\documentclass[a4paper,fleqn,sort&compress]{cas-sc}
\usepackage[authoryear,longnamesfirst]{natbib}
\usepackage{hyperref}
\usepackage{amssymb}
\usepackage{amsmath}
\allowdisplaybreaks[4]
\usepackage{amssymb}
\usepackage{epstopdf}
\usepackage{booktabs}
\usepackage[subfigure]{tocloft}
\usepackage{subfigure}
\usepackage{graphicx}
\usepackage{indentfirst}
\usepackage{float}
\usepackage{amsthm}
\newtheorem{theorem}{Theorem}
\theoremstyle{thmstyletwo}%
\newtheorem{lemma}{Lemma}
\newdefinition{rmk}{Remark}
\newenvironment{pf}{\noindent{\bf Proof.}}{\scriptsize{$\quad \square$}  \par}

\newdefinition{Definition}{Definition}
\newtheorem{Proposition}{Proposition}

\newtheorem{Property}{Property}
\begin{document}
\shortauthors{Zhang et~al.}
\let\printorcid\relax
\let\WriteBookmarks\relax
\def\floatpagepagefraction{1}
\def\textpagefraction{.001}
\title [mode = title]{Prescribed-time boundary control of second-order hyperbolic PDEs modeled flexible string systems via backstepping design}                      
\tnotemark[1]


\author[1]{Chuan Zhang}

\author[1]{He Yang}

\author[1]{Fei Wang}	

\affiliation[1]{organization={School of Mathematical Science},
	addressline={Qufu Normal University},
	city={Qufu},
	postcode={273165},
	state={Shandong},
	country={China}}

\author[2,3]{Tuo Zhou}
\cormark[1]
\cortext[2]{Corresponding author}
\ead{tzhou@sdtbu.edu.cn}
\affiliation[2]{organization={School of Information and Electronic Engineering},
	addressline={Shandong Technology and Business University},
	city={Yantai},
	postcode={264005},
	state={Shandong},
	country={China}}
\affiliation[3]{organization={Research Institute of Information Fusion},
	addressline={Naval Aviation University},
	city={Yantai},
	postcode={264005},
	state={Shandong},
	country={China}}

\begin{abstract}
	This paper presents a boundary control scheme for prescribed-time (PT) stable of flexible string systems via backstepping method, and the dynamics of such systems modeled by Hamilton's principle is described as second-order hyperbolic partial differential equations (PDEs). Initially, to construct a boundary controller with PT stabilization capacity, a PT stable hyperbolic PDEs system with time-varying coefficient is chosen as the target system, and a corresponding Volterra integral transform with time-varying kernel function is considered. Then, to identify the boundary controller, the well-posedness of kernel equation is derived by means of successive approximation and mathematical induction, and the upper bound of kernel function is estimated. Furthermore, the inverse transform is proved with the help of a similar process for kernel function. Subsequently, the PT stability of closed-loop system is proved by PT stability of target system and reversible integral transform. Finally, the simulation results demonstrate the effectiveness of our scheme.
\end{abstract}




\maketitle
\section{Introduction}\label{sec1}
In practical engineering environments, nearly all mechanical structures exhibit flexibility, necessitating advanced control strategies for flexible systems. Such systems are ubiquitous in industrial automation, marine engineering, and aerospace applications, where their dynamic behavior significantly impacts operational efficiency and safety \cite{liu2017modeling,zhu2019enhanced,he2022trajectory}. To model these systems, three primary methodologies are employed: Hamilton's principle for energy-based formulations \cite{zhang2005pde}, Newton-Euler equations for force-moment analysis \cite{scaglioni2017closed}, and Lagrange equations for generalized coordinate dynamics \cite{kang2002dynamic}. Unlike rigid-body systems, flexible structures exhibit infinite-dimensional dynamics, typically modeled as distributed parameter systems (DPSs) governed by partial differential equations (PDEs).\\
\indent Control strategies for DPS generally follow two paradigms. The first involves discretizing the infinite-dimensional model into a finite-dimensional approximation for controller synthesis. However, this approach risks spillover instability due to neglected high-frequency modes, potentially degrading closed-loop performance \cite{hagen2003spillover}. The second paradigm directly addresses the infinite-dimensional nature of the system, employing methodologies such as distributed control \cite{mo2024distributed,wang2023distributed}, boundary control \cite{chentouf2020exponential}, or pointwise control \cite{song2023improved}. Boundary control, which actuates the system at its spatial endpoints, offers practical advantages over distributed control by reducing sensor/actuator density, lowering implementation costs, and simplifying operation. For instance, He et al. \cite{he2020dynamical} utilized boundary control to suppress torsional-flexural vibrations in flexible aircraft wings while ensuring angular positioning via Lyapunov analysis. Similarly, Zhao et al. \cite{zhao2019boundary} developed a boundary anti-disturbance controller for nonlinear flexible strings under unknown perturbations. Numerous scholars have investigated boundary control problems for PDEs, with controllability results have been established in \cite{lions1988exact,glowinski1995exact}. Recent advances include event-triggered boundary control for delayed reaction-diffusion systems \cite{koudohode2024event}, which guarantees exponential stability while minimizing actuation updates. Among these methods, backstepping method proposed by Krstic stands out for its systematic handling of nonlinearities through Volterra integral transformations \cite{krstic2008backstepping,de2024backstepping,krstic2008boundary}.\\
\indent The backstepping methodology stabilizes unstable DPS by constructing a Volterra integral transformation that maps the original system to a target stable system. Critical to this approach is solving the associated kernel equations to derive boundary control laws, with rigorous analysis of solution existence and well-posedness. Recent applications include sliding-mode backstepping for parabolic DPS with mixed boundary conditions \cite{bao2023sliding}. Hou et al. \cite{hou2024state} further extended this framework to fractional-order DPS with time delays. While these studies achieve asymptotic or exponential stabilization, boundary control for PT stabilization that convergence occurs exactly at a prescribed time for hyperbolic PDEs such as flexible string is still unexplored, which motivates this research. \\
\indent For some time-constrained applications, such as missile guidance \cite{jain2024comparison}, multi-agent formation control \cite{meurer2011finite} and spacecraft rendezvous \cite{cortes2006finite}, demand convergence within strict or predefined time. Although the finite-time (FT) control exhibits a quicker convergence speed compared to asymptotic convergence, its convergence time depends on the initial conditions \cite{ghaderi2024output}, complicating practical deployment. Fixed-time (FxT) control decouples convergence time from initial states but requires conservative upper bounds tied to system parameters \cite{tao2022fixed}. In contrast, PT control ensures exact stabilization at a designer-specified time, independent of initial conditions or tuning parameters \cite{ning2023novel,wei2022prescribed,ye2022prescribed}. Wei et al. \cite{wei2022prescribed} investigated the problem of PT stabilization for a class of semilinear parabolic systems subject to spatiotemporal-varying disturbance via distributed control. Xiao et al. \cite{xiao2024prescribed} investigated the PT constrained tracking control problem for a class of 2×2 hyperbolic PDE systems actuated by nonlinear ODE dynamics. In contrast, this work addresses a second-order hyperbolic PDE system and proposes a PT boundary control scheme that is more direct and computationally efficient. By avoiding the delays inherent in coupled or distributed control strategies, the method remains inherently compatible with wave propagation constraints. Stabilization requires a prescribed time $T > T_{min}$, where  $T_{min}$ corresponds to the wave propagation round-trip time. Our proposed PT boundary control explicitly respects this limit while enabling exact convergence at $T = T_{min}+\epsilon$, where $\epsilon>0$.\\
\indent This paper addresses the open issue of PT stabilization for hyperbolic PDEs modeling flexible string systems via boundary backstepping control. The key contributions are:
\begin{itemize}
	\item While existing backstepping designs have achieved PT stabilization for parabolic PDEs systems \cite{wei2024prescribed} and specific hyperbolic PDEs system \cite{xiao2024prescribed}, critical gaps remain in addressing flexible string systems modeled by second-order hyperbolic PDE. This paper bridges this gap by proposing a PT boundary control scheme for string systems modeled by second-order hyperbolic PDEs. Given the challenges in obtaining analytical solutions for the time-varying kernel function, we adopt a numerical solution approach, effectively addressing the dynamic characteristics of the system. The designed boundary control demonstrates good robustness in stabilizing the string system within a specified time frame. 
	\item Different from the distributed control strategy utilized in \cite{mo2024distributed,wang2023distributed}, the boundary control proposed in this paper is physically more realistic when applied to string systems due to the actuation and sensor are all nonintrusive, and it does not require to construct additional state observers.
	\item By analyzing the well-posedness of kernel equation and inverse kernel equation, we establish the equivalence of PT stability between the original system and the target system. Additionly, the boundary controller is explicitly related to the kernel function, and an algorithm is given to obtain the numerical expression of kernel function.
\end{itemize}
\indent \textbf{Notation.} In this paper, $\mathbb{R}$ denotes the set of real numbers, $\mathbb{R}^{+}$ expresses the non-negative real numbers. $\mathbb{N}$ is the set of non-negative integers, $\mathbb{C}^{\infty }$ denotes infinitely continuously differentiable. $\frac{\partial ^{l} }{\partial t^{l} }$ represents the $l$-order time partial derivative, $\frac{d^{l} }{dt^{l} }$ denotes $l$-order time derivative. $L^{2}\left [0,1 \right ] $ denotes the set of functions on the interval $\left [ 0,1 \right ]$ satisfying $\int_{0}^{1}p^{2}\left ( x,t \right )dx< \infty$.  $\left \| p\left ( x,t \right )  \right \|_{L^{2} }$ denotes the norm of $p\left ( x,t \right ) \in L^{2}\left [0,1 \right ]$ and $\left \| p\left ( x,t \right )  \right \|_{L^{2} }=\left ( \int_{0}^{1} p^{2} \left ( x,t \right ) dx \right )^{\frac{1}{2} }$.
\section{Mathematical Modeling and Preliminary}\label{sec2}
\begin{figure}[!h]
	\centering
	\includegraphics[width=5.3cm]{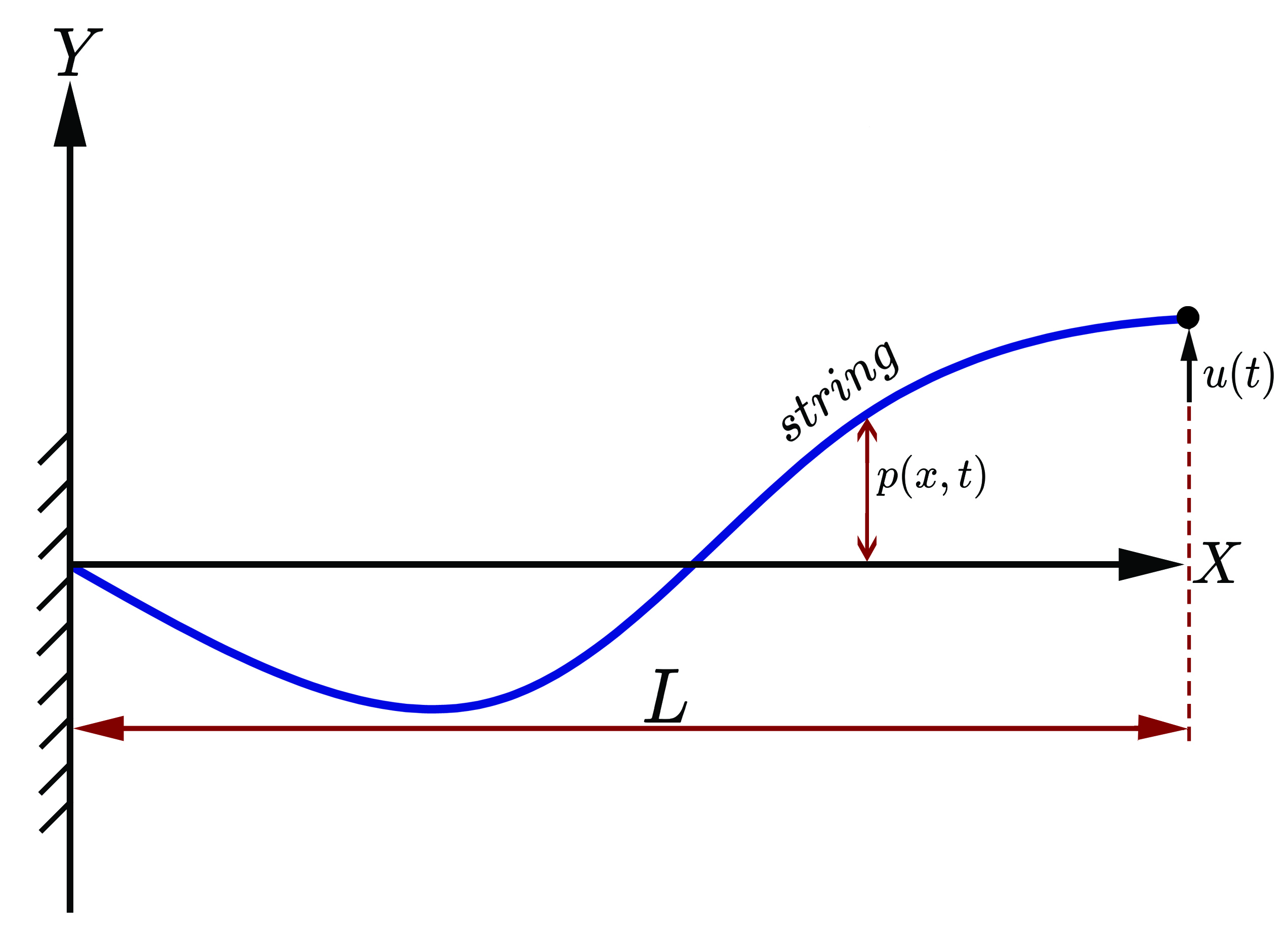}
	\caption{The flexible string systems}
	\label{fig.1}
\end{figure}

This paper takes the flexible string system on the two-dimensional plane as the research object, and Fig.\ref{fig.1} shows the simplified model of this system. To simplify the description, the symbols $\left ( \cdot  \right )=\left ( \cdot  \right )\left ( t \right ) $ and $\left ( \cdot  \right )=\left ( \cdot  \right )\left (x, t \right ) $ are employed throughout the paper. $\frac{dV\left ( t \right ) }{dt}$, $\frac{\partial p\left ( x,t\right ) }{\partial x}$, $\frac{\partial p^{2} \left ( x,t\right ) }{\partial x^{2} }$, $\frac{\partial p \left ( x,t\right ) }{\partial t }$ and $\frac{\partial p^{2}  \left ( x,t\right ) }{\partial t^{2}  }$ are denoted by $\dot{V}\left ( t \right )$, $p_{x} \left ( x,t \right ) $, $p_{xx} \left ( x,t \right ) $, $p_{t} \left ( x,t \right ) $ and $p_{tt} \left ( x,t \right ) $, respectively. \\
\indent The left boundary of the string is assumed to be fixed. $p\left ( 1,t \right )$, $p_{t} \left ( 1,t \right )$ and $p_{tt} \left ( 1,t \right )$ are the displacement, velocity and acceleration of the tip payload respectively. $u\left ( t \right )$ is the boundary control input; Let $\rho _{0}$ be the uniform mass per unit length of the string; $M$ denotes the mass of the payload at the right boundary of the string; $T_{f}$ is a constant tension; $x\in \left [ 0,1 \right ]$ is the spatial independent variable, and $t\in \left [ 0,T\right ) $ is the time independent variable.\\
\indent The kinetic energy of the string system is expressed as
	\begin{equation*}
		E_{k} \left ( t \right )=\frac{1}{2}M\left [ p_{t}\left ( 1,t \right )   \right ] ^{2}+\frac{1}{2}\rho _{0}\int_{0}^{1}\left [ p_{t}\left ( x,t \right )\right ]^{2}dx,        
	\end{equation*}
the potential energy due to the constant tension $T_{f}$ is denoted as
	\begin{equation*}
		E_{p} \left ( t \right )=\frac{1}{2}T_{f} \int_{0}^{1}\left [ p_{x}\left ( x,t \right )\right ]^{2}dx,        
	\end{equation*}
the total virtual work done by the system is
	\begin{equation*}
		\delta W\left ( t \right ) =u\left ( t \right )\delta p\left ( 1,t \right ).      
	\end{equation*}
\indent We can derive the equations of motion describing the flexible string from Hamilton's principle as follows:
	\begin{equation*}
		\int_{t_{1}}^{t_{2}}\delta\left[E_{k}\left(t\right)-E_{p}\left(t\right)+W\left(t\right)\right ]dt=0, 
	\end{equation*}
where $t_{1}$, $t_{2}$ are two time variables; $\delta$ denotes the variational operator; $E_{k}\left ( t \right )$ and $E_{p}\left ( t \right )$ represent kinetic and potential energy, respectively. $W\left ( t \right )$ represents the virtual power done by the non-conservative forces acting on the system. This yields the following hyperbolic PDEs system:
	\begin{subequations}\label{shi1}
		\begin{align}
			\rho_{0}\cdot p_{tt}\left(x,t\right)&=T_{f}\cdot p_{xx}\left(x,t\right),\label{shi1a} \\   
			p\left ( 0,t \right )&=0, \label{shi1b}  \\
			T_{f}\cdot  p_{x}\left ( 1,t \right )+M\cdot p_{tt}\left ( 1,t \right )&=u\left ( t \right ),   \label{shi1c}   \\
			p\left ( x,0 \right ) &=p_{0}\left ( x \right ), \label{shi1d} 
		\end{align}
	\end{subequations}
where (\ref{shi1b}) and (\ref{shi1c}) are the boundary conditions, (\ref{shi1d}) is the initial condition.
\begin{Proposition}(Minimal Controllability Time):\label{Pro.1}
	For the flexible string system governed by hyperbolic PDE, there exists a minimum controllable time $T_{min}$ determined by the wave propagation speed. The minimal time required for boundary controllability of the string system is given by: 
	\begin{equation*}
		T_{min}=2*\frac{L}{c} ,
	\end{equation*}
	where $L=1$ is the string length, and $c=\sqrt{T_{f}/\rho_{0}}$ is the wave speed derived from the PDE (\ref{shi1a}), $T_{min}$ denotes the time of one round trip of the wave.
\end{Proposition}
\begin{Definition}\label{def.1}
	The controlled system (\ref{shi1}) is PT stable if there exists a boundary controller $u\left ( t \right )$ such that for any initial value $p_{0}\left ( x \right ) \in L^{2}\left [ 0,1 \right ]$, the solution of system (\ref{shi1}) satisfies 
	\begin{equation*}
		\left \| p\left ( \cdot ,t \right )  \right \|_{L^{2}\left[ 0,1 \right ]  }\to 0,\quad t\to T
	\end{equation*}
		where $x\in \left [ 0,1 \right ]$, $t\in \left [ 0,T\right )$, $T>T_{min}$ is prescribed time.
	\end{Definition}
	\begin{Definition}\label{def.2}
		The first class of $\lambda$-order modified Bessel function $I_{\lambda }\left ( \cdot  \right )$ is defined as 
			\begin{equation*}
				I_{\lambda }\left ( f \right )=\sum_{n=0}^{\infty }\frac{1}{n!\Gamma \left ( n+\lambda +1 \right ) }\left ( \frac{f}{2}  \right ) ^{\lambda +2n}=\sum_{n=0}^{\infty }\frac{1}{n!\left (  n+\lambda \right )!  }\left ( \frac{f}{2}  \right ) ^{\lambda +2n}.
			\end{equation*}
	\end{Definition}
	\begin{Property}\label{A1}
		\cite{yu2024adaptive}
		If the kinetic energy of the system (\ref{shi1}) is bounded, then $p_{t}\left ( x,t \right )$ and $p_{xt}\left ( x,t \right )$ are bounded for all $ \left ( x,t \right )\in \left [ 0,1 \right ]\times \left [ 0,T \right )$.
	\end{Property}
	\begin{Property}\label{A2}
		\cite{yu2024adaptive}
		If the potential energy of the system (\ref{shi1}) is bounded, then $p_{x}\left ( x,t \right )$ and $p_{xx}\left ( x,t \right )$ are bounded for all $\left ( x,t \right )\in \left [ 0,1 \right ]\times \left [ 0,T \right )$.
	\end{Property}
	\begin{lemma}\label{L1}
		\cite{meurer2009tracking}
		For non-negative integers $n\in \mathbb{N}$, $j\in \mathbb{N}$, $l\in \mathbb{N}$, $j\le l$, the following equation is satisfied:
			\begin{equation*}
					\sum_{j=0}^{l} \binom{l}{j}\left ( l-j+1 \right )!\left ( j+n \right )!=\sum_{j=0}^{l}\binom{l}{j}\left ( l-j+n \right )!\left ( j+1 \right )!=\frac{\left ( l+n+2 \right )! }{\left ( n+1 \right ) \left ( n+2 \right ) }, 
			\end{equation*}
		where $l\in \mathbb{N}$, $n\in \mathbb{N}$.
	\end{lemma}
	\begin{lemma}\label{L2}
		\cite{meurer2009tracking}
		For any real numbers $\xi > 0$ and $\eta > 0$, the following equation holds:
			\begin{equation*}
				\begin{aligned}
					\int_{0}^{\eta } \int_{0}^{\tau } \left ( s\tau  \right )^{n-1}\left ( s+\tau  \right ) dsd\tau &=\frac{\eta ^{2n+1} }{n\left ( n+1 \right ) },\quad \eta > 0\\
					\int_{\eta }^{\xi } \int_{0}^{\eta } \left ( s\tau  \right )^{n-1}\left ( s-\tau  \right ) dsd\tau& =\frac{\left ( \xi \eta  \right )^{n} \left ( \xi -\eta  \right )  }{n\left ( n+1 \right ) }, \quad 0\le \eta \le \xi 
				\end{aligned}
			\end{equation*}
		where $n\in \mathbb{N}$, $n\ge 1$.
	\end{lemma}
	
	\begin{lemma}\label{L3}
		\cite{ge2009boundary,zhao2018output}
		Let $\psi _{1}\left ( x,t \right ) $, $\psi _{2}\left ( x,t \right ) \in \mathbb{R} $, then the following inequality holds:
			\begin{equation*}
				\begin{aligned}
					\psi _{1}\psi _{2} &\le \lvert \psi _{1}\psi _{2} \rvert\le \psi _{1}^{2} + \psi _{2}^{2}, \\
					\lvert \psi _{1}\psi _{2} \rvert &=\lvert \left ( \frac{1}{\sqrt{\delta } } \psi _{1}  \right ) \left ( \sqrt{\delta } \psi _{2}  \right ) \rvert \le \frac{1}{\delta }\psi _{1}^{2}+\delta \psi _{2}^{2},\quad \delta > 0.
				\end{aligned}
			\end{equation*} 
	\end{lemma}	
	\begin{lemma}\label{L4}
		\cite{rahn2001mechatronic}
		For $\hbar\left ( x,t \right )\in\mathcal{H}^1[0,L]$, if the condition $\hbar\left ( 0,t \right ) =0$, $\forall t\in \left [ 0,T\right )$ is satisfied, then it satisfies:
			\begin{equation*}
				\hbar ^{2}\le L\cdot\int_{0}^{L}\left [\hbar {'}\right ]^{2}dx,\quad \int_{0}^{L} \hbar ^{2}dx\le L^{2} \cdot \int_{0}^{L}\hbar _{x}^{2}dx.
			\end{equation*}
	\end{lemma}
	\section{PT boundary controller design}\label{sec3}
	\subsection{Selection and PT stability of target system}\label{sec3.1}
	Unlike traditional Volterra transforms with time-invariant kernels, our design incorporates a time-varying kernel $k\left ( x,y,t \right )$ scaled by $\mu\left (t \right )$, enabling PT convergence by coupling the system dynamics with a time-varying coefficient:
		\begin{equation}\label{shi2}
			\upsilon \left ( x,t \right ) =p\left ( x,t \right )-\int_{0}^{x}k\left ( x,y,t \right )p\left ( y,t \right ) dy,
		\end{equation}
with this transform (\ref{shi2}), we convert the original system (\ref{shi1}) into the following target system:
	\begin{subequations}\label{shi3}
		\begin{align}
			\rho _{0}\cdot  \upsilon _{tt} \left ( x,t \right )&=T_{f}\cdot  \upsilon _{xx} \left ( x,t \right )-\mu \left ( t \right )\upsilon \left ( x,t \right ),    \label{shi3a}  \\
			\upsilon\left ( 0,t \right )&=0,  \label{shi3b}   \\
			T_{f}\cdot \upsilon_{x} \left (1,t \right )&=-M \cdot \upsilon_{tt}  \left ( 1,t \right ) , \label{shi3c}  \\
			\upsilon\left (x,0 \right )&=\upsilon _{0} \left ( x \right ),  \label{shi3d}
		\end{align}
	\end{subequations}
where $\mu \left ( t \right ) \in \mathcal{C^{\infty } } \left [ 0,T \right ) $. To stabilize the target system within the prescribed time, $\mu \left ( t \right )$ is designed as:
	\begin{equation}\label{shi4}
		\mu \left ( t \right )=\frac{\mu _{0}^{2}T^{2}}{\left ( T-t \right ) ^{2} }, \quad \mu _{0}> 0, \quad t\in\left [ 0,T\right ).
	\end{equation}		
\indent From the above equation, $\mu \left ( t \right )$ is monotonically increasing, $\mu \left ( 0\right )=\mu _{0}^{2}$ and $\lim\limits_{t \to T}\mu \left ( t \right ) =+\infty$. Furthermore, taking the $l\in \mathbb{N}$-order derivative of the time-varying function $\mu \left ( t \right ) $ with respect to time $t$ yields $\frac{d^{l} }{dt^{l} }\mu \left ( t \right ) =\mu ^{\frac{l}{2}+1 }\left ( t \right )\left ( l+1 \right )!\frac{1}{\left ( \mu _{0}T  \right )^{l}}$. 
\begin{rmk}
	This is different from the case in \cite{wei2022prescribed} where the parameter design is associated with the diffusion coefficient in the parabolic PDEs, in this paper, our hyperbolic PDEs formulation only requires $\mu _{0}>0$. Our work complements \cite{wei2022prescribed} by extending the PT stabilization framework to hyperbolic PDEs.
\end{rmk}	
\begin{theorem}\label{Th1}
	Consider the target system (\ref{shi3}) with the time-varying coefficient $\mu \left ( t \right )$ defined in (\ref{shi4}), for any prescribed time $T>T_{min}$, where $T_{min}$ is the minimal controllable time, and any initial condition $\upsilon_{0}\left(x\right )\in L^{2}\left[0,1\right]$, the $L^{2}$-norm of $\upsilon\left(x,t\right )$ converges to zero as $t\to T$:
		\begin{equation*}
				\left\|\upsilon\left(\cdot ,t\right)\right\|_{L^{2}}\le\sqrt{\frac{2}{\sigma _{2} T_{f} } \left \{  \varsigma \left(t\right)\left(V\left(0\right)+\varepsilon\cdot\int_{0}^{t}e^{\frac{2\lambda _{1} }{\sigma _{3}}\int_{0}^{t}\mu\left(s\right)ds}ds\right) \right \}},
		\end{equation*}
	where $\lim\limits_{t \to T}\varsigma \left ( t \right )=0 $. 
\end{theorem}
\begin{pf}
	To prove that the closed-loop system (\ref{shi1}) is PT stabilization under the PT boundary controller, it is initially necessary to prove the PT stabilization of the target system. \\
	\indent Consider the Lyapunov candidate function
		\begin{equation}\label{shi5}
			V\left ( t \right )=V_{1} \left ( t \right )+V_{2} \left ( t \right ), 
		\end{equation}
	in which
		\begin{align}
			V_{1}\left(t\right)&=\frac{\beta}{2}\rho_{0}\int_{0}^{1}\upsilon_{t}^{2}\left(x,t\right)dx+\frac{\beta}{2}T_{f} \int_{0}^{1}\upsilon_{x}^{2}\left(x,t\right)dx, \label{shi6} \\
			V_{2}\left(t\right)&=\alpha\rho_{0}\int_{0}^{1}x\upsilon_{x}\left(x,t\right)\upsilon_{t}\left(x,t\right)dx,  \label{shi7} 
		\end{align}
	where $t\in \left [ 0,T\right )$, $\alpha$ and $\beta$ are two positive weighting constants.   
	\begin{lemma}\label{L5}
		The upper and lower bounds of the Lyapunov candidate function given by (\ref{shi5}) as
			\begin{equation*}
				0\le \sigma _{2}V_{1}\left ( t \right )\le V\left ( t \right )\le \sigma _{3}V_{1}\left ( t \right ),
			\end{equation*}
		where $V_{1}\left ( t \right )$ is given in (\ref{shi6}), $\sigma_{2}$ and $\sigma_{3}$ are constants defined as
			\begin{equation*}
					\sigma _{2}=1-\frac{2\alpha\rho_{0}}{\min \left ( \beta \rho_{_{0}},\beta T_{f} \right ) }> 0, \quad
					\sigma _{3}=1+\frac{2\alpha\rho_{0}}{\min \left ( \beta \rho_{_{0}},\beta T_{f} \right ) }> 1.
			\end{equation*}
		\begin{pf}
			By Lemma \ref{L3}, we know (\ref{shi7}) can be transformed into
		\begin{equation*}
					\lvert V_{2}\left(t\right)\rvert\le \alpha \rho _{0}\int_{0}^{1}\left(\upsilon _{t}^{2}\left(x,t\right) +\upsilon _{x}^{2}\left(x,t\right)\right)dx\le \sigma_{1}V_{1}\left(t\right),
		\end{equation*}
			where $\sigma_{1} =\frac{2\alpha \rho_{0}}{\min\left(\beta\rho_{0},\beta T_{f}\right)} $, and we can further get 
		\begin{equation*}
					-\sigma_{1}V_{1}\left(t\right)\le V_{2}\left(t\right)\le \sigma_{1}V_{1}\left(t\right ),
		\end{equation*}
			given that $\alpha$ is a small positive weighting constant satisfying $0<\alpha<\frac{\min\left(\beta\rho _{0},\beta T_{f}\right)}{2\rho_{0}}$, we can get
				\begin{equation*}
						\sigma _{2}=1-\sigma _{1} =1-\frac{2\alpha\rho_{0}}{\min \left ( \beta \rho_{_{0}},\beta T_{f} \right ) }> 0,\quad
						\sigma _{3}=1+\sigma _{1} =1+\frac{2\alpha\rho_{0}}{\min \left ( \beta \rho_{_{0}},\beta T_{f} \right ) }> 1,
				\end{equation*}
			therefore, we further get 
				\begin{equation*}
					0\le\left(1-\sigma_{1}\right)V_{1}\le V_{1}\left(t\right)+V_{2}\left(t\right)=V\left(t\right)\le\left(1+\sigma_{1}\right)V_{1} \left (t \right ),
				\end{equation*}
			which means 
				\begin{equation*}
					0\le \sigma _{2} V_{1}\le V\left(t\right)\le \sigma _{3} V_{1}\left(t\right ).
				\end{equation*}
		\end{pf}
	\end{lemma}
	\indent Deriving (\ref{shi6}) and (\ref{shi7}) along the trajectory of the target system (\ref{shi3}) with respect to time $t$
		\begin{equation*}
			\begin{aligned}
				\dot{V_{1}}\left(t\right)&=\beta T_{f}\int_{0}^{1}\upsilon_{x}\left(x,t\right)\upsilon_{xt}\left(x,t\right )dx+\beta \rho _{0} \int_{0}^{1}\upsilon_{t}\left(x,t\right)\upsilon_{tt}\left(x,t\right)dx \notag \\
				&=\beta T_{f} \upsilon_{x}\left(1,t\right)\upsilon_{t}\left(1,t\right)-\beta T_{f} \int_{0}^{1}\upsilon_{t}\left(x,t\right )\upsilon_{xx}\left(x,t\right)dx+\beta T_{f} \int_{0}^{1}\upsilon_{t}\left(x,t\right )\upsilon_{xx}\left(x,t\right)dx\notag\\
				&\quad-\beta\mu\left(t\right)\int_{0}^{1}\upsilon_{t}\left(x,t\right)\upsilon\left(x,t\right)dx   \notag  \\
				&=\beta T_{f} \upsilon_{x}\left(1,t\right)\upsilon_{t}\left(1,t\right)-\beta\mu\left(t\right)\int_{0}^{1}\upsilon_{t}\left(x,t\right )\upsilon\left(x,t\right)dx ,   \\
			\end{aligned}
		\end{equation*}
		\begin{equation*}\label{shi38}
			\begin{aligned}
				\dot{V_{2}}\left(t\right)&=\alpha\rho _{0} \int_{0}^{1}\left(x\upsilon_{xt}\left(x,t\right)\upsilon_{t}\left(x,t\right)+x\upsilon_{x}\left(x,t\right)\upsilon_{tt}\left(x,t\right)\right)dx     \notag   \\
				&=\alpha\rho _{0} \int_{0}^{1}\left(x\upsilon_{xt}\left(x,t\right)\upsilon_{t}\left(x,t\right)+x\upsilon_{x}\left(x,t\right)\cdot\left ( \frac{T_{f} }{\rho _{0} } \upsilon_{xx}\left(x,t\right)-\frac{1}{\rho _{0} } \mu \left ( t \right ) \upsilon \left ( x,t \right ) \right ) \right )dx    \notag     \\
				&=\alpha\rho_{0}\int_{0}^{1}x\upsilon_{xt}\left(x,t\right)\upsilon_{t}\left(x,t\right)dx-\alpha\int_{0}^{1}x\upsilon_{x}\left(x,t\right)\mu\left(t\right)\upsilon\left(x,t\right)dx  +\alpha T_{f} \int_{0}^{1}x\upsilon_{x}\left(x,t\right)\upsilon_{xx}\left(x,t\right)dx,           
			\end{aligned}
		\end{equation*}
	by Lemma \ref{L3} and the method of integration by parts yields:
		\begin{align}
			\dot{V_{1}}\left(t\right)&\le\beta T_{f} \upsilon_{x}\left(1,t\right)\upsilon_{t}\left(1,t\right)+\frac{\beta\mu\left(t\right)}{\delta_{1}}\int_{0}^{1}\upsilon_{t}^{2}\left(x,t\right)dx +\beta\mu\left(t\right)\delta_{1}\int_{0}^{1}\upsilon^{2} \left(x,t\right)dx \label{shi8}\\
			\dot{V_{2}}\left(t\right)&\le\frac{\alpha T_{f} }{2}\upsilon_{x}^{2}\left(1,t\right)-\frac{\alpha T_{f} }{2}\int_{0}^{1}\upsilon_{x}^{2}\left(x,t\right)dx+\frac{\alpha\mu\left(t\right)}{\delta_{2}}\int_{0}^{1}\upsilon^{2}\left(x,t\right)dx \notag  \\
			&\quad+\alpha\mu\left(t\right)\cdot\delta_{2}\int_{0}^{1}\upsilon_{x}^{2}\left(x,t\right)dx+\frac{\alpha\rho _{0} }{2}\upsilon_{t}^{2}\left(1,t\right)-\frac{\alpha\rho _{0} }{2}\int_{0}^{1}\upsilon_{t}^{2}\left(x,t\right)dx\label{shi9}
		\end{align}
	\indent Applying the partial derivative of (\ref{shi5}) relative to time $t$, and substituting (\ref{shi8}) and (\ref{shi9}) yields:
		\begin{equation*}
			\begin{aligned}
				\dot{V}\left(t\right)&=\dot{V}_{1}\left(t\right)+\dot{V_{2}}\left(t\right)   \\
				&\le\beta T_{f} \upsilon_{x}\left(1,t\right)\upsilon_{t}\left(1,t\right)+\frac{\beta\mu\left(t\right)}{\delta_{1}}\int_{0}^{1}\upsilon_{t}^{2}\left(x,t\right)dx+\beta\mu\left(t\right)\delta_{1}\int_{0}^{1}\upsilon^{2}\left(x,t\right)dx+\frac{\alpha T_{f}}{2}\upsilon_{x}^{2}\left(1,t\right)\\
				&\quad-\frac{\alpha T_{f}}{2}\int_{0}^{1}\upsilon_{x}^{2}\left(x,t\right)dx+\frac{\alpha\mu\left(t\right)}{\delta_{2}}\int_{0}^{1}\upsilon^{2}\left(x,t\right)dx+\alpha\mu\left(t\right)\delta_{2}\int_{0}^{1}\upsilon_{x}^{2}\left(x,t\right)dx\\
				&\quad+\frac{\alpha\rho_{0}}{2}\upsilon_{t}^{2}\left(1,t\right)-\frac{\alpha\rho_{0} }{2}\int_{0}^{1}\upsilon_{t}^{2}\left(x,t\right)dx  \\
				&\le-\left(\frac{\alpha\rho _{0}  }{2\mu\left(t\right)}-\frac{\beta}{\delta_{1}}\right)\mu\left(t\right)\int_{0}^{1}\upsilon_{t}^{2}\left(x,t\right)dx-\left(\frac{\alpha T_{f} }{2\mu\left(t\right)}-\alpha\delta_{2}-\beta\delta_{1}-\frac{\alpha }{\delta _{1} }\right)\mu\left(t\right)\int_{0}^{1}\upsilon_{x}^{2}\left(x,t\right)dx\\
				&\quad+\frac{\delta_{3} }{2}\left(\upsilon_{x}\left(1,t\right)+\upsilon_{t}\left(1,t\right)\right)^{2}\\
				&\le-\lambda _{1} \mu\left(t\right)V_{1}\left(t\right)+\varepsilon,   \\
			\end{aligned}
		\end{equation*}	
	where
		\begin{equation*}
			\left\{\begin{matrix}
				\varphi _{1}=\frac{\alpha \rho _{0} }{2\mu \left ( t \right ) }-\frac{\beta }{\delta _{1} }>0,  \\
				\varphi _{2}=\frac{\alpha T_{f} }{2\mu \left ( t \right ) }-\alpha \delta _{2}-\beta \delta _{1}  -\frac{\alpha }{\delta _{1} }>0 ,\\       
				\lambda _{1}=\min \left ( \frac{2\varphi _{1}}{\beta  \rho _{0} } ,\frac{2\varphi _{2}}{\beta  T_{f} }  \right )> 0,\\
				\delta _{3}=\max \left ( \alpha  T_{f} ,\alpha \rho _{0} ,\beta  T_{f}  \right )> 0,\\		
				\varepsilon=\frac{\delta_{3}}{2}\left[\upsilon_{x}\left(1,t\right)+\upsilon_{t}\left(1,t\right)\right]^{2}>0.
			\end{matrix}\right.       
		\end{equation*}
	\indent Through the Sobolev embedding theorem and the definition of $V_{1}$, we know that both $V_{x}$ and $V_{t}$ are governed by $V$. Given that $V$ is exponentially dominated by $\mu$, the growth rate of $\varepsilon$ remains strictly bounded by the decay term, thereby ensuring system stability.\\
	\indent from Lemma \ref{L5}, we can further deduce:
		\begin{equation}\label{shi10}
			\dot{V}\left(t\right)\le-\frac{2\lambda _{1} }{\sigma _{3}}\mu\left(t\right )V\left(t\right)+\varepsilon.
		\end{equation}
	\indent Integrating both sides of (\ref{shi10}) with respect to time $t$ from 0 to $t$, one can compute:
		\begin{equation}\label{shi11}
			\begin{aligned}
				V\left(t\right)\le e^{-\frac{2\lambda _{1} }{\sigma _{3}}\int_{0}^{t}\mu\left(s\right)ds}V\left(0\right)+\varepsilon\cdot e^{-\frac{2\lambda _{1} }{\sigma _{3}}\int_{0}^{t}\mu\left(s\right)ds}\cdot\int_{0}^{t}e^{\frac{2\lambda _{1} }{\sigma _{3}}\int_{0}^{t}\mu\left(s\right)ds}ds,
			\end{aligned}
		\end{equation}
	let $\varsigma \left(t\right)=e^{-\frac{2\lambda _{1} }{\sigma _{3}}\int_{0}^{t}\mu\left(s\right)ds}$ and use $\mu\left(t\right)=\frac{\mu_{0}^{2}T^{2}}{\left(T-t\right)^{2}}$, $\mu_{0}>0$, we can get
	\begin{equation*}
			\varsigma \left(t\right)=e^{\frac{2\lambda _{1} }{\sigma _{3}}\mu_{0}^{2}T}e^{-\frac{2\lambda _{1}}{\sigma _{3}}\mu_{0}T\sqrt{\mu\left(t\right)}}=e^{\frac{2\lambda _{1}}{\sigma _{3}}\mu_{0}^{2}T}e^{-\frac{2\lambda _{1}}{\sigma _{3}}\mu_{0}^{2}T\left(\frac{T}{T-t}\right)},
	\end{equation*}
	where $t\in\left[0,T\right)$, this function is monotonically decreasing and satisfies $\varsigma \left(0\right)=1$ and $ \displaystyle\lim_{t\to T} \varsigma\left(t\right)=0$, which implies that $V\left(t\right)$ is bounded. Using Lemma \ref{L4}, Lemma \ref{L5} and (\ref{shi6}), for $t\in\left[0,T\right)$, there is
		\begin{equation}\label{shi12}
			\begin{aligned}
				\frac{T_{f} }{2}\left\|\upsilon\left(x,t\right)\right\|_{L^{2}}^{2}=\frac{T_{f}}{2}\int_{0}^{1}\upsilon^{2}\left(x,t\right)dx\le\frac{T_{f}}{2}\int_{0}^{1}\upsilon_{x}^{2}\left(x,t\right)dx\le V_{1}\left(t\right)\le\frac{1}{\sigma _{2}}V\left(t\right)\in \mathcal{L}_{\infty},
			\end{aligned}
		\end{equation}
	from this we can deduce
		\begin{equation}\label{shi13}
			\begin{aligned}
				\left\|\upsilon\left(\cdot,t\right)\right\|_{L^{2}}	\le\sqrt{\frac{2}{\sigma _{2}T_{f} }V\left(t\right)}\le\sqrt{\frac{2}{\sigma _{2} T_{f} } \left \{  \varsigma \left(t\right)\left[V\left(0\right)+\varepsilon\int_{0}^{t}e^{\frac{2\lambda _{1} }{\sigma _{3}}\int_{0}^{t}\mu\left(s\right)ds}ds\right] \right \}},
			\end{aligned}
		\end{equation}
	we can deduce that $\displaystyle\lim_{t\to T}\left \| \upsilon \left (\cdot,t \right )  \right \|_{L^{2} } \to 0 $ as $t\to T$, from the Definition \ref{def.1}, the target system (\ref{shi3}) satisfies the PT stability in the time $T>T_{min}$.  
\end{pf}
\begin{rmk}
	The PT stability analysis presented in this work is also applicable to the stability analysis of Lyapunov-based controller designs.
\end{rmk}
\subsection{Determination and well-posedness of kernel equations}\label{sec3.2}
\indent To derive the PT boundary controller, it is necessary to obtain the kernel functions in (\ref{shi2}). Following the backstepping methodology \cite{lei2020stabilization}, we first compute the temporal derivatives of the transformation (\ref{shi2}):
	\begin{equation}\label{shi14}
		\begin{aligned}
			\upsilon _{t}\left ( x,t \right )=p_{t}\left ( x,t \right )-\int_{0}^{x}k_{t}\left ( x,y,t \right )p\left ( y,t \right )dy-\int_{0}^{x}k\left ( x,y,t \right )p_{t}\left ( y,t \right ) dy.
		\end{aligned}
	\end{equation}
	\begin{equation*}
			\upsilon _{tt}\left ( x,t \right )=p_{tt}\left ( x,t \right )-\int_{0}^{x}k_{tt}\left ( x,y,t \right )p\left ( y,t \right )dy-2\int_{0}^{x} k_{t}\left ( x,y,t \right ) p_{t} \left ( y,t \right ) dy- \int_{0}^{x}k\left ( x,y,t \right )p_{tt}\left ( y,t \right ) dy.   
	\end{equation*}	
Substituting (\ref{shi1a}) into the above equation and applying integration by parts yields:
	\begin{equation}\label{shi15}
		\begin{aligned}		
			\upsilon _{tt}\left ( x,t \right )
			&=\frac{T_{f} }{\rho _{0} } p_{xx} \left ( x,t \right )-\int_{0}^{x}k_{tt}\left ( x,y,t \right )p\left ( y,t \right )dy-2\int_{0}^{x} k_{t}\left ( x,y,t \right ) p_{t} \left ( y,t \right ) dy\\
			&\quad-\frac{T_{f} }{\rho _{0} } k\left ( x,x,t \right )p_{x}\left ( x,t \right )+\frac{T_{f} }{\rho _{0} } k\left ( x,0,t \right )p_{x}\left ( 0,t \right )+\frac{T_{f} }{\rho _{0} } p\left ( x,t \right )k_{y}\left ( x,x,t \right )\\
			&\quad-\frac{T_{f} }{\rho _{0} } p\left ( 0,t \right )k_{y}\left ( x,0,t \right )-\frac{T_{f} }{\rho _{0} } \int_{0}^{x}p\left ( y,t \right )k_{yy}\left ( x,y,t \right )dy,  
		\end{aligned}
	\end{equation}
	\indent Next, spatial derivatives of (\ref{shi2}) are computed using Leibniz's differential rule \cite{liu2003boundary}:
		\begin{align}
			\upsilon _{x}\left ( x,t \right )&=p_{x}\left ( x,t \right )-k\left ( x,x,t \right )p\left ( x,t \right )-\int_{0}^{x}k_{x}\left ( x,y,t \right )p\left ( y,t \right )dy,\notag\\
			\upsilon _{xx}\left ( x,t \right )&=p_{xx}\left ( x,t \right )-\frac{d}{dx} k\left ( x,x,t \right )p\left ( x,t \right )-k\left ( x,x,t \right )p_{x}\left ( x,t \right )\notag\\
			&\quad-k_{x}\left ( x,x,t \right )p \left ( x,t \right )-\int_{0}^{x}k_{xx}\left ( x,y,t \right )p\left ( y,t \right )dy. \label{shi16}
		\end{align}
	where $k_{x}\left ( x,x,t \right )=k_{x}\left ( x,y,t \right )\mid _{y=x}$, $\frac{d}{dx} k\left ( x,x,t \right )=k_{x}\left ( x,y,t \right )+k_{y}\left ( x,y,t \right )$.\\
	\indent Let $x = 1$ and combine (\ref{shi1c}) and (\ref{shi3c}) to get:
		\begin{align*}
			T_{f}\cdot\upsilon_{x}\left(1,t\right)+M\cdot\upsilon_{tt}\left(1,t\right)&=T_{f}\cdot p_{x}\left(1,t\right)-T_{f}\cdot k\left(1,1,t\right)p\left(1,t\right)-T_{f}\int_{0}^{1}k_{x}\left ( 1,y,t \right )p\left ( y,t \right ) dy\\
			&\quad+M\cdot p_{tt}\left ( 1,t \right )-\frac{MT_{f} }{\rho _{0}}k\left ( 1,1,t \right ) p_{x}\left ( 1,t \right )+\frac{MT_{f} }{\rho _{0}}k_{y}\left ( 1,1,t \right )p\left ( 1,t \right )\\
			&\quad-\frac{MT_{f} }{\rho _{0}}\int_{0}^{1}k_{yy}\left ( 1,y,t \right )p\left ( y,t \right )dy=0, 
		\end{align*}
	to ensure that the original system (\ref{shi1}) is transformed into the target system (\ref{shi3}) by the integral transform (\ref{shi2}), we have to choose the following controller:
		\begin{equation}\label{shi17}
			\begin{aligned}
				u\left ( t \right ) &=T_{f}\cdot k\left(1,1,t\right)p\left(1,t\right)+T_{f}\int_{0}^{1}k_{x}\left ( 1,y,t \right )p\left ( y,t \right ) dy+\frac{MT_{f} }{\rho _{0}}k\left ( 1,1,t \right ) p_{x}\left ( 1,t \right )\\
				&\quad-\frac{MT_{f} }{\rho _{0}}k_{y}\left ( 1,1,t \right )p\left ( 1,t \right )+\frac{MT_{f} }{\rho _{0}}\int_{0}^{1}k_{yy}\left ( 1,y,t \right )p\left ( y,t \right )dy.
			\end{aligned}
		\end{equation}
	\indent Finally, substituting (\ref{shi2}), (\ref{shi15}) and (\ref{shi16}) into the target system (\ref{shi3}) yields
		\begin{equation}\label{shi18}
			\begin{aligned}
			&\rho _{0}\cdot  \upsilon_{tt}  \left ( x,t \right )-T_{f}\cdot  \upsilon_{xx}  \left ( x,t \right )+\mu \left ( t \right )\upsilon \left ( x,t \right )\\&=T_{f}\cdot  p_{xx}\left ( x,t \right ) -\rho _{0}\int_{0}^{x}k_{tt}\left ( x,y,t \right )p\left ( y,t \right )dy-2\rho _{0} \int_{0}^{x}k_{t} \left ( x,y,t \right ) p_{t} \left ( y,t \right ) dy-T_{f}\cdot   k\left ( x,x,t \right )p_{x}\left ( x,t \right )\\
			&\quad+T_{f}\cdot  k\left ( x,0,t \right )p_{x}\left ( 0,t \right )+T_{f}\cdot p\left ( x,t \right )k_{y}\left ( x,x,t \right ) -T_{f}\cdot p\left ( 0,t \right )k_{y}\left ( x,0,t \right )    \\
			&\quad-T_{f}\cdot \int_{0}^{x}p\left ( y,t \right )k_{yy}\left ( x,y,t \right )dy-T_{f}\cdot p_{xx}\left ( x,t \right )+T_{f}\cdot \frac{d}{dx}k\left ( x,x,t \right )p\left ( x,t \right ) \\
			&\quad+T_{f}\cdot k\left ( x,x,t \right )p_{x}\left ( x,t \right )+T_{f}\cdot k_{x}\left ( x,x,t \right )p\left ( x,t \right )+T_{f}\int_{0}^{x}k_{xx}\left ( x,y,t \right )p\left ( y,t \right )dy\\
			&\quad+\mu \left ( t \right )\left ( p\left ( x,t \right ) -\int_{0}^{x}k\left ( x,y,t \right )p\left ( y,t \right ) dy\right)      \\
			&=\left (T_{f}\cdot  k_{y} \left ( x,x,t \right ) +T_{f}\cdot\frac{d}{dx}k\left ( x,x,t \right ) +T_{f}\cdot k_{x}\left ( x,x,t \right ) +\mu \left ( t \right )   \right )\cdot p\left ( x,t \right )\\
			&\quad+\int_{0}^{x}\left ( -\rho _{0} \cdot k_{tt}\left ( x,y,t \right )-T_{f} \cdot k_{yy}\left ( x,y,t \right )+T_{f} \cdot k_{xx}\left ( x,y,t \right )-\mu \left ( t \right )k\left ( x,y,t \right )\right )p\left ( y,t \right )dy\\&\quad-2\rho _{0} \int_{0}^{x}k_{t} \left ( x,y,t \right ) p_{t} \left ( y,t \right ) dy+T_{f} \cdot k\left ( x,0,t \right )p_{x}\left ( 0,t \right )-T_{f}\cdot  p\left ( 0,t \right )k_{y}\left ( x,0,t \right )=0. 
			\end{aligned}
		\end{equation}
	\indent Substituting the boundary condition (\ref{shi1b}) into (\ref{shi18}) yields the kernel equations with time-varying coefficient:
		\begin{subequations}\label{shi19}
			\begin{align}
				\rho _{0} k_{tt}\left ( x,y,t \right )&=-T_{f} k_{yy}\left ( x,y,t \right )+T_{f} k_{xx}\left ( x,y,t \right )-\mu \left ( t \right )k\left ( x,y,t \right ),  \\ 
				\frac{d}{dx}k\left ( x,x,t \right )&=-\frac{1}{2T_{f} }\mu \left ( t \right ),    \\      
				k\left ( x,0,t \right )&=0,  
			\end{align}
		\end{subequations}
	where the domain of definition of kernel functions $k\left ( x,y,t \right )$ is $\mathbb{D}=\left \{ \left ( x,y,t \right )\in\mathbb{R }^{2}\times \left [ 0,T \right ):0\le y\le x\le 1\right \}$, the initial condition for the kernel equations (\ref{shi19}) is 
	\begin{equation*}
		\int_{0}^{x}k\left ( x,y,0\right )p_{0}\left ( y \right ) dy=p_{0}\left ( x \right )-\upsilon _{0}\left ( x \right ). 
	\end{equation*}
	\begin{rmk}
		This paper solves the problem of well-posedness of the kernel equations in the case of no analytic solution, and the proposed method is also applicable to the proof of well-posedness of kernel equations in \cite{espitia2019boundary}, which is more generalized.
	\end{rmk}
	\begin{rmk}
		Since the kernel equations contain a function $\mu \left ( t \right )\in \mathcal{C}^{\infty } \left [ 0,T\right )$ with property $\lim\limits_{t \to T}\mu \left ( t \right )=+\infty$, it is important not only to prove the well-posedness of kernel equations (\ref{shi19}), but also to obtain an explicit expression for the upper bound estimate of the time-varying kernel functions with respect to the time correlation.
	\end{rmk}
	\begin{rmk}
		Note that in contrast to the general kernel equations \cite{su2017boundary}, the kernel equations in this paper are hyperbolic distributed parameter systems with time-varying coefficient, and it is difficult to obtain their analytic solution. The uniqueness and existence of kernel function solution is proved in the Appendix. In Algorithm 1, the numerical solution of the kernel functions is given based on the finite difference method, and the controller is obtained using the numerical solution.
	\end{rmk}
	\begin{theorem}\label{Th2}
		Considering (\ref{shi1}) and the target system (\ref{shi3}), the function $\mu \left ( t \right )\in \mathcal{C}^{\infty } \left [ 0,T\right )$ is designed as (\ref{shi4}), then there is a unique solution to the kernel equations (\ref{shi19}) on the domain of definition $\mathbb{D}$, and the upper bound of the kernel functions $k\left ( x,y,t \right )$ is estimated as
			\begin{equation}\label{shi20}
					\lvert k\left ( x,y,t \right ) \rvert\le  \frac{y\left ( 6+\mu _{0}^{2} T^{2}\right)}{T_{f} }\cdot \mu \left ( t \right )\cdot \left ( \frac{e\mu \left ( t \right )\left ( 6+\mu _{0}^{2} T^{2}\right )\left ( x^{2} -y^{2}  \right ) }{T_{f}\mu _{0} ^{2}T^{2}}\right )^{-\frac{1}{2} }I_{1}\left ( z \right ),
			\end{equation}
		where $\mu _{0}>0$, $T$ is prescribed time. $I_{1}\left ( f \right )$ denotes the first class of first-order modified Bessel function.
	\end{theorem}
		\begin{table}[htbp]
		\begin{tabular*}{\hsize}{@{}@{\extracolsep{\fill}}l@{}}
			\toprule[1.5pt]  %
			\textbf{Algorithm 1}: The steps of the numerical solution algorithm for the kernel function\\
			\midrule
			\textbf{Step 1}: Give the initialization parameters $T$, $\mu _{0}$, $\rho _{0}$, $T_{f}$, $M$ based on Theorem \ref{Th4};\\
			\textbf{Step 2}: Discretize $x$, $y$ and $t$, determine the space step $hx$, $hy$ and the time step $ht$, and generate the sequence of\\ discrete points $i$, $j$, $n$, and compute the number of discrete points xnumber, ynumber, tnumber;\\
			\textbf{Step 3}: Initialize the matrix used to store $k\left ( x,y,t \right )$ and its derivatives;\\
			\textbf{Step 4}: Discrete $k_{xx}$, $k_{yy}$, $k_{tt}$, construct difference equations, use for loops to compute $p\left ( x,0 \right )$, $\mu \left ( t \right )$, $k \left ( 1,1,t \right )$ and\\ $k\left ( x,y,0 \right )$; use nested loops to compute the values of $p\left ( x,y\right ) $, $k\left ( x,y,t\right )$ at each discrete point;\\
			\textbf{Step 5}: The value obtained in step 4 is used in a nested loop for $u\left ( t \right )$. The outer loop traverses the time step for\\ accumulating the integral approximation, and the inner loop traverses the space step for approximating the integral\\ calculation;\\
			\bottomrule
		\end{tabular*}	
	\end{table}
	\indent \indent The proof of Theorem \ref{Th2} is given in the Appendix. Observe that the time-varying gain kernel grows unboundedly near the end time $T$ , but the boundary controller is bounded on $t\in \left [ 0,T\right )$, as will be shown later in the PT stabilization of the closed-loop system.\\
	\indent To demonstrate that the closed-loop system (\ref{shi1}) and the target system (\ref{shi3}) can be converted to each other by the integral transform (\ref{shi2}), it is necessary to prove that the integral transform (\ref{shi2}) is invertible.\\
	\indent The inverse transformation of the integral transform (\ref{shi2}) is given by
	\begin{equation}\label{shi21}
		p\left ( x,t \right )=\upsilon \left ( x,t \right )+\int_{0}^{x} r\left ( x,y,t \right )\upsilon \left ( y,t \right )dy, 
	\end{equation}
	where the domain of definition of $r\left ( x,y,t \right )$ is $\mathbb{D}$. The integral transform (\ref{shi2}) and the inverse transform (\ref{shi21}) satisfy the below relation
	\begin{equation}\label{22}
		r\left ( x,y,t \right )=k\left ( x,y,t \right )+\int_{y}^{x}k\left ( x,\gamma ,t \right ) r\left ( \gamma ,y,t \right )d\gamma, 
	\end{equation}
	the inverse transformation (\ref{shi21}) converts the target system (\ref{shi3}) into the original system (\ref{shi1}) with the boundary controller (\ref{shi17}). Substituting the integral transform (\ref{shi21}) into the original system (\ref{shi1}), and then utilizing the target system (\ref{shi3}), the method of integration by parts and Leibniz's differential rule, the kernel equations can be obtained:
		\begin{subequations}\label{shi23}
			\begin{align}
				\rho _{0}\cdot  r_{tt}\left(x,y,t\right)&=-T_{f}\cdot  r_{yy}\left(x,y,t\right)+T_{f}\cdot  r_{xx}\left(x,y,t\right)+\mu\left(t\right)r\left ( x,y,t \right ) ,  \label{shi32a}  \\             
				\frac{d}{dx} r\left(x,x,t\right)&=-\frac{1}{2T_{f} }\mu \left ( t \right ), \quad x\in \left [ 0,1 \right ],\label{shi32b} \\
				r\left(x,0,t\right)&=0 , \quad x\in \left [ 0,1 \right ],\label{shi32c} 
			\end{align}
		\end{subequations}
	where $\left(x,y,t\right)\in \mathbb{D}$.  \\
	\indent The kernel equations (\ref{shi23}) are similar to (\ref{shi19}), then by Theorem \ref{Th2}, there is also a unique solution to the kernel equations (\ref{shi23}), and the inverse kernel functions $r\left ( x,y,t \right )$ have the same properties as the kernel functions $k\left ( x,y,t \right )$. It is evident that we can obtain the below Theorem:
	\begin{theorem}\label{Th3}
		Considering the original system (\ref{shi1}) and the target system (\ref{shi3}), $\mu\left(t\right)\in \mathcal{C}^{\infty}\left [ 0,T\right )$ is taken as (\ref{shi4}), then there is a unique solution of the kernel equations (\ref{shi23}) on $\mathbb{D}$, which is second-order continuously differentiable with respect to the space $x$, $y$ and infinitely continuously differentiable with respect to the time $t$. Moreover, the upper bound of the kernel functions $r\left ( x,y,t \right )$ is estimated as
			\begin{equation}\label{shi24}
					\lvert r\left(x,y,t\right)\rvert\le \frac{y\left ( 6+\mu _{0}^{2} T^{2}\right)}{T_{f} }\cdot \mu \left ( t \right )\cdot \left ( \frac{e\mu \left ( t \right )\left ( 6+\mu _{0}^{2} T^{2}\right )\left ( x^{2} -y^{2}  \right ) }{T_{f}\mu _{0} ^{2}T^{2}}\right )^{-\frac{1}{2} }I_{1}\left ( z \right ),
			\end{equation}
		where $0\le y\le x\le 1$, $t\in \left [ 0,T\right ) $, and $T$ is prescribed time.
	\end{theorem}
	\begin{rmk}
		The proof procedure for the inverse kernel function is similar to the proof of Theorem \ref{Th2}.
	\end{rmk}
	\subsection{PT stability of closed-loop system}\label{sec3.3}
	The target system (\ref{shi3}) is adapted on a compact subset $\left [ T_{min},T^{*} \right]$, $T_{min}<T^{*}<T$. Therefore, from the well-posedness of the target system and the kernel equations as well as the bounded and inverse integral transform (\ref{shi2}) on $\left [T_{min},T^{*} \right]$, the closed-loop system (\ref{shi1}) based on the PT boundary controller (\ref{shi17}) is well-posedness for $ t\in \left [ 0,T\right )$.
	\begin{theorem}\label{Th4}
		If there exists a monotonically increasing  function $\mu \left ( t \right )$ in (\ref{shi4}) such that for all $t\in \left [ 0,T\right )$,
			\begin{equation*}
					\left\|p\left(\cdot,t\right)\right\|_{L^{2}}\le\left(1+\mathfrak{M}_{r}\left(t\right)\right)\cdot\sqrt{\frac{2}{\sigma _{2} T_{f} }\left [ \varsigma \left(t\right)\left(V\left(0\right)+\varepsilon\int_{0}^{t}e^{\frac{2\lambda _{1} }{\sigma _{3}}\int_{0}^{t}\mu\left(s\right)ds}ds\right) \right ] },
			\end{equation*}
		where $\varsigma \left ( t \right ) \to 0$ as $t\to T$, then for any initial condition $p_{0}\left(x\right )\in L^{2}\left[0,1\right]$, the closed-loop system (\ref{shi1}) with the boundary controller (\ref{shi17}) is PT stable in the prescribed time.
\end{theorem}
\begin{pf}
	Since the kernel functions $k\left(x,y,t\right)$ and $r\left(x,y,t\right)$ are continuous on $\mathbb{D}$, we can obtain
		\begin{equation*}
			\begin{aligned}
				\left\|k\left(x,y,t\right)\right\|_{\infty}=\underset{0\le y\le x\le 1}{\sup}\lvert k\left(x,y,t\right)\rvert&\le \mathfrak{M}_{k}\left(t\right ),  \quad t\in\left[0,T\right).     \\ \left\|r\left(x,y,t\right)\right\|_{\infty}=\underset{0\le y\le x\le 1}{\sup}\lvert r\left(x,y,t\right)\rvert&\le \mathfrak{M}_{r}\left(t\right ), \quad t\in\left[0,T\right).
			\end{aligned}
		\end{equation*}
	\indent In order to determine the explicit expressions for $\mathfrak{M}_{k}\left(t\right )$ and $\mathfrak{M}_{r}\left(t\right )$, using the power series property
		\begin{equation*} \displaystyle\sum_{n=0}^{\infty}\Phi_{n}^{m}\le\left(\sum_{n=0}^{\infty}\Phi_{n}\right)^{m}, \quad m\ge1, \quad \Phi_{n}\ge0, \quad n\in\mathbb{N} 
		\end{equation*}
	and the first class of first-order modified Bessel function
		\begin{equation*}
			\begin{aligned}
			 I_{1}\left(\sqrt{\frac{\mu\left(t\right)\cdot e\cdot \left(6+\mu _{0}^{2}T^{2}\right)\left(x^{2}-y^{2}\right) }{T_{f}\mu _{0} ^{2} T^{2}  } } \right)=\sum_{n=0}^{\infty}\frac{\left( \displaystyle \frac{\sqrt{\mu\left(t\right)\cdot e\cdot \left(6+\mu_{0}^{2}T^{2}\right)\left(x^{2}-y^{2}\right)/T_{f}\mu _{0} ^{2} T^{2} } }{2}\right)^{2n+1}}{n!\left(n+1\right)!}, 
			\end{aligned}
		\end{equation*}
	which results in
		\begin{equation}\label{shi25}
			\begin{aligned}
				&\frac{1}{2} \left ( \frac{e\cdot \mu \left ( t \right )\left ( 6+\mu _{0}^{2} T^{2}\right) \left ( x^{2}-y^{2}   \right )  }{T_{f}\mu _{0}^{2}T^{2} }  \right ) ^{-\frac{1}{2} }\cdot \sum_{n=0}^{\infty}\frac{\left( \displaystyle \frac{\sqrt{\mu\left(t\right)\cdot e\cdot \left(6+\mu_{0}^{2}T^{2}\right)\left(x^{2}-y^{2}\right)/T_{f}\mu _{0} ^{2} T^{2} } }{2}\right)^{2n+1}}{n!\left(n+1\right)!}\\
				&\le\sum_{n=0}^{\infty}\frac{1}{\left ( n! \right ) ^{2} }\left(\frac{\sqrt{\mu\left(t\right)\left(6+\mu_{0}^{2}T^{2}\right)\left(x^{2}-y^{2}\right)/T_{f}\mu _{0} ^{2} T^{2} } }{2}\right)^{2n}  \\
				&=\left (\sum_{n=0}^{\infty }\left ( \frac{\sqrt{e\mu \left ( t \right )\left ( 6+\mu _{0} ^{2} T ^{2}\right ) \left ( x^{2}-y^{2} \right )  /T_{f} \mu _{0} ^{2} T^{2}} }{2}  \right )^{n}/n!\right )^{2}= e^{\frac{\sqrt{e\cdot \mu\left(t\right)\left(6+\mu _{0}^{2}T^{2}\right)\left(x^{2}-y^{2}\right)/T_{f} } }{\mu_{0}T}}.
			\end{aligned}
		\end{equation}
	\indent Then, by (\ref{shi25}) it is possible to obtain
		\begin{equation}\label{shi26}
				\lvert \frac{y\left ( 6+\mu _{0}^{2}T^{2}\right ) }{T_{f} } \mu \left ( t \right ) \frac{I_{1}\left ( \sqrt{\mu \left ( t \right )e\left(6+\mu_{0} ^{2}T^{2}\right)\left(x^{2}-y^{2}\right)/T_{f}}/\mu _{0}T \right)}{\sqrt{\mu \left ( t \right ) e \left ( 6+\mu_{0}^{2} T^{2}   \right )  \left ( x^{2} -y^{2}  \right ) /T_{f} } /\mu _{0}T } \rvert\le C_{k}\cdot \mu \left ( t \right )\cdot e^{\frac{\sqrt{e\cdot \mu \left ( t \right ) \left ( 6+\mu _{0}^{2}T  ^{2} \right ) /T_{f} } }{\mu _{0}T } },
		\end{equation}
	where $C_{k}>0$ is a constant. From the estimate (\ref{shi20}) of the upper bound on the kernel function $k\left ( x,y,t \right )$, an explicit expression for $\mathfrak{M}_{k}\left(t\right)$ can be obtained as
		\begin{equation}\label{shi27}
			\lvert k\left(x,y,t\right)\rvert\le \mathfrak{M}_{k}\left(t\right)=C_{k}\cdot \mu \left ( t \right ) e^{\frac{\sqrt{e\cdot \mu \left ( t \right )\left ( 6+\mu _{0}^{2}T^{2}\right ) /T_{f}  } }{\mu _{0} T} },
		\end{equation}
	similarly, the estimate (\ref{shi24}) of the upper bound on the inverse kernel function $r\left ( x,y,t \right )$ leads to the explicit expression for $\mathfrak{M}_{r}\left ( t \right )$
		\begin{equation}\label{shi28}
			\lvert r\left(x,y,t\right)\rvert\le \mathfrak{M}_{r}\left(t\right)=C_{k}\cdot \mu \left ( t \right ) e^{\frac{\sqrt{e\cdot \mu \left ( t \right )\left ( 6+\mu _{0}^{2}T^{2}\right ) /T_{f}  } }{\mu _{0} T} }. 
		\end{equation}
	\indent Using the inverse integral transforms (\ref{shi2}) and (\ref{shi21}), and the Cauchy-Schwarz inequality, combined with estimates of the upper bounds on the kernel functions $k\left ( x,y,t \right )$ and $r\left ( x,y,t \right )$, it can be obtained that
		\begin{align}
			\left\|p\left(\cdot,t\right)\right\|_{L^{2}}&\le\left(1+\mathfrak{M}_{r}\left(t\right)\right)\cdot\left\|\upsilon\left(\cdot,t\right)\right\|_{L^{2} },\label{shi29} \\
			\left\|\upsilon\left(\cdot,t\right)\right\|_{L^{2}}&\le\left(1+\mathfrak{M}_{k}\left(t\right)\right)\cdot\left\|p\left(\cdot,t\right)\right\|_{L^{2} }, \label{shi30}
		\end{align}
	substituting (\ref{shi13}) into (\ref{shi29}) yields
		\begin{equation}\label{shi31}
			\begin{aligned}
				\left\|p\left(\cdot,t\right)\right\|_{L^{2}}&\le\left(1+\mathfrak{M}_{r}\left(t\right)\right)\cdot\left\|\upsilon\left(\cdot,t\right)\right\|_{L^{2}}   \\
				&\le\left(1+\mathfrak{M}_{r}\left(t\right)\right)\cdot\sqrt{\frac{2}{\sigma _{2} T_{f} }\left [ \varsigma \left(t\right)\left(V\left(0\right)+\varepsilon\int_{0}^{t}e^{\frac{2\lambda _{1} }{\sigma _{3}}\int_{0}^{t}\mu\left(s\right)ds}ds\right) \right ] }, 
			\end{aligned}
		\end{equation}
	using (\ref{shi28}) and (\ref{shi31}), for $t\in\left[0,T\right)$, it follows that
		\begin{equation*}
				\left\|p\left(\cdot,t\right)\right\|_{L^{2}}
				\le\left(1+C_{k}\cdot\mu \left ( t \right )\cdot e^{\frac{\sqrt{e\cdot \mu \left ( t \right )\left ( 6+\mu _{0}^{2}T^{2}\right ) /T_{f}  } }{\mu _{0}T} } \right) \cdot\sqrt{\frac{2}{\sigma _{2}T_{f}  } \left [\varsigma \left(t\right)\left(V\left(0\right)+\varepsilon\int_{0}^{t}e^{\frac{2\lambda _{1} }{\sigma _{3}}\int_{0}^{t}\mu\left(s\right)ds}ds\right) \right ] },
		\end{equation*}
	then substituting the function $\mu\left(t\right)=\frac{\mu_{0}^{2}T^{2}}{\left (T-t\right )^{2} }$ into the above equation we have
		\begin{equation}\label{shi32}
				\left\|p\left(\cdot,t\right)\right\|_{L^{2}}
				\le\left(1+C_{k}\cdot\frac{\mu_{0}^{2}T^{2}}{\left (T-t\right )^{2} }\cdot e^{\frac{\sqrt{e\left ( 6+\mu _{0}^{2}T^{2}\right ) /T_{f}  } }{T-t} } \right) \cdot\sqrt{\frac{2}{\sigma _{2}T_{f}  } \left [\varsigma \left(t\right)\left(V\left(0\right)+\varepsilon\int_{0}^{t}e^{\frac{2\lambda _{1} }{\sigma _{3}}\int_{0}^{t}\mu\left(s\right)ds}ds\right) \right ] },
		\end{equation}
	we know that $\displaystyle\lim_{t\to T}\varsigma \left(t\right)=0$ when $t\to T$, so we have $\displaystyle\lim_{t\to T}\left\|p\left(\cdot,t\right)\right\|_{L^{2}}=0$. According to the Definition \ref{def.1}, the closed-loop system (\ref{shi1}) is PT stabilization in the prescribed time $T > T_{min}$.    
	\begin{rmk}
		From (\ref{shi12}), we get that $V_{1}\left ( t \right ) $ is bounded on $t\in \left [ 0,T\right )$, and thus $\upsilon _{t}\left ( x,t \right )$ and $\upsilon _{x}\left ( x,t \right )$ are bounded on $\forall \left ( x,t \right )\in \left [ 0,1 \right ]\times \left [ 0,T \right )$, and further $p_{x}\left ( x,t \right )$, $p_{t}\left ( x,t \right )$ are bounded on $\forall \left ( x,t \right )\in \left [ 0,1 \right ]\times \left [ 0,T \right )$, then we get that $E_{k}\left (t\right)$ and $E_{p}\left (t\right)$ are bounded. Using Properties \ref{A1} and \ref{A2}, we can obtain that $p_{xx}\left (x,t\right)$ and $p_{xt}\left (x,t\right)$ are also bounded. Finally, we use (\ref{shi1a}) to see that $p_{tt}\left (x,t\right)$ is also bounded. Since $p\left (x,t\right)$, $p_{x} \left (x,t\right)$, and $k\left (x,y,t\right)$ are all bounded, we can conclude that the boundary controller is also bounded.
	\end{rmk} 
\end{pf}
\begin{rmk}
	The proposed PT control design inherits a fundamental trade-off between exact convergence and noise robustness, which is intrinsic to time-varying gain approach. As demonstrated for parabolic PDEs \cite{wei2024prescribed} and linear systems \cite{aldana2023inherent}, the time-varying gain $\mu\left (t\right)$ inevitably decreases noise robustness. Some workarounds have been proposed in \cite{song2017time} and \cite{holloway2019prescribed}, \cite{song2017time} proposed another solution for switching off the algorithm when the error trajectory enters the desired error deadzone; \cite{holloway2019prescribed} proposed to switch off the algorithm at  time $t_{stop} $ before the prescribed time $T$, thereby keeping the time-varying gain uniformly bounded. Due to the nature of these algorithms, unperturbed trajectories with large initial conditions will enter such deadzone arbitrarily close to $T$ while having arbitrarily large is time-varying gain. Furthermore, in the presence of additional arbitrarily small noise, the trajectory may not enter the deadzone at all, thus lacking robustness despite the workaround. In practical applications, we should adjust the design parameters carefully to obtain suitable transient performance and control effects. In future, we will further explore how to reduce the sensitivity to noise while maintaining PT convergence.
\end{rmk}
\section{Numerical simulations}\label{sec4}
\indent To validate the effectiveness of the proposed boundary control scheme, numerical simulations are conducted for both the open-loop and closed-loop systems. The flexible string system is discretized using the finite difference method with the following parameters: mass per unit length $\rho _{0}=1kg/m$, tip payload mass $M=1kg$, constant tension $T_{f}=45N$, time-varying coefficient $\mu_{0}=5$, prescribed time $T=3s>T_{min}$; the initial condition is $p_{0}\left ( x \right )=-\frac{1}{2} x\left ( x-1 \right )$.
\begin{figure}[!h]
	\centering
	\includegraphics[width=7.0cm]{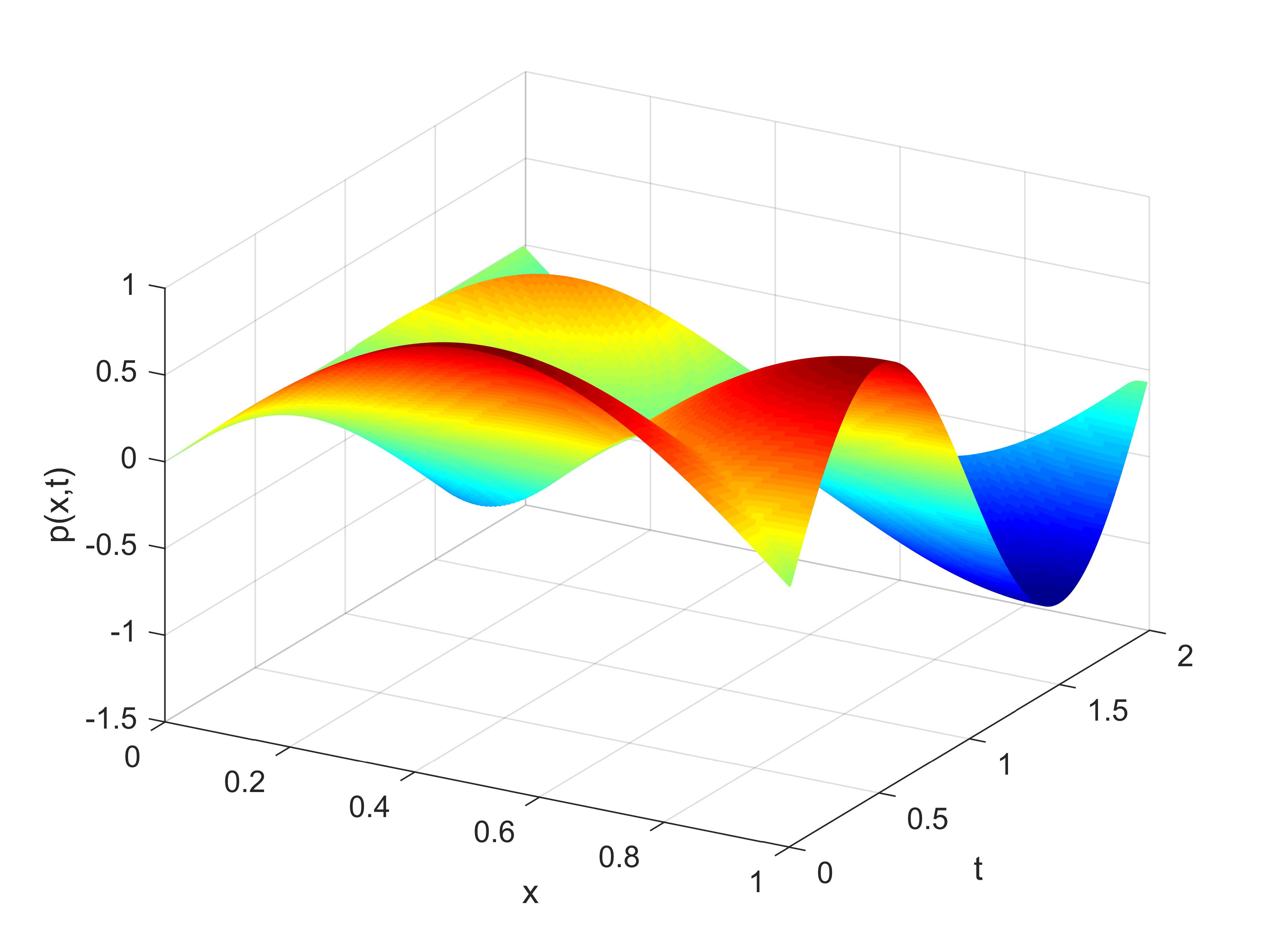}
	\caption{The evolution of the open-loop system state}
	\label{fig.2}
\end{figure}\par
\begin{figure}[!h]
	\centering
	\includegraphics[width=7.0cm]{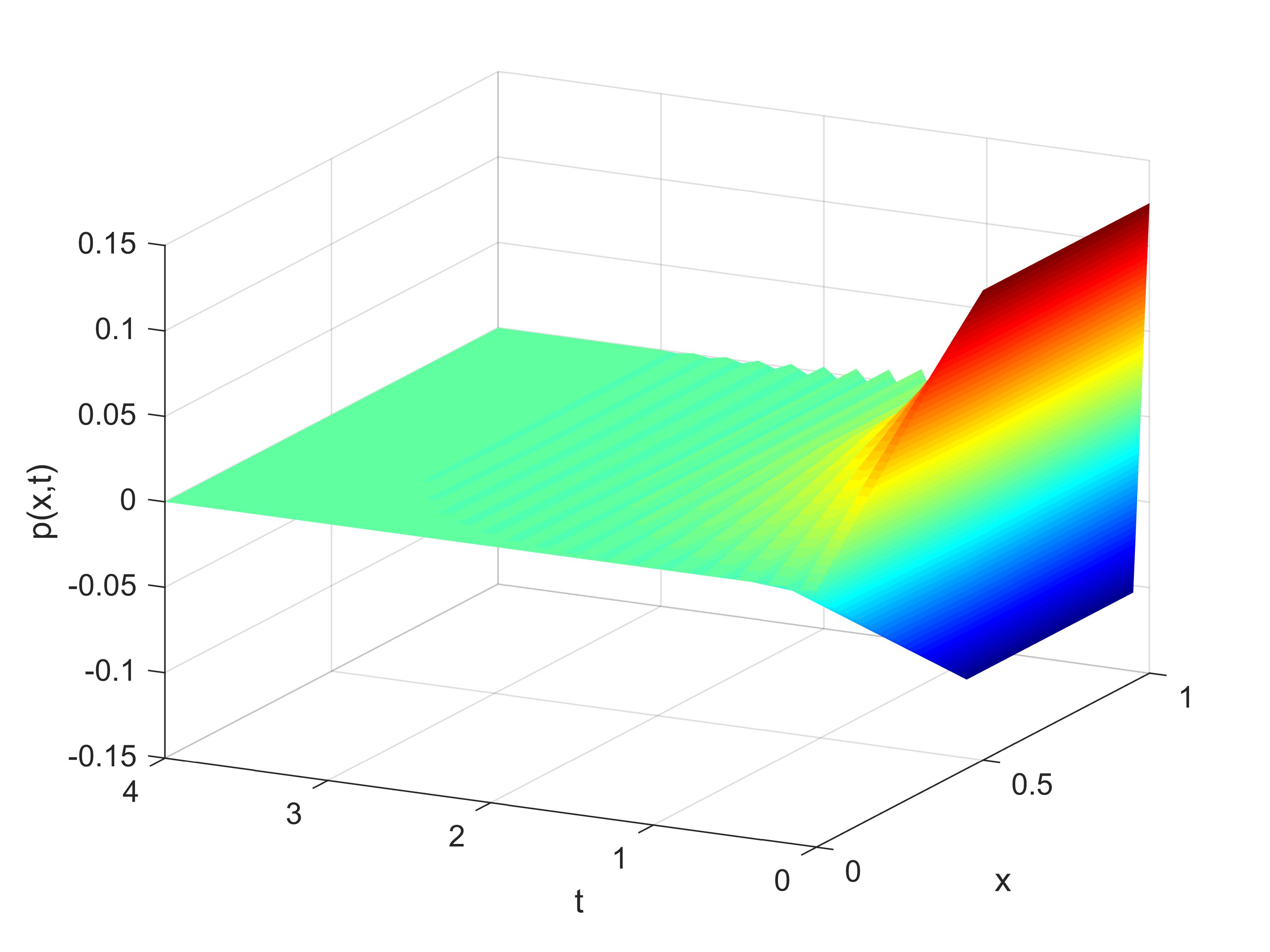}
	\caption{The evolution of the closed-loop system state}
	\label{fig.3}
\end{figure}\par
\indent Fig.\ref{fig.2} illustrates the spatiotemporal evolution of the open-loop system state. The system exhibits unstable oscillatory behavior due to the lack of control input, confirming the necessity of stabilization. Fig.\ref{fig.3} demonstrates the system state evolution with the proposed boundary controller (\ref{shi17}). The result shows that the closed-loop system achieves stabilization within the prescribed time, validating the PT convergence property of the controller. 
\begin{figure}[!h]
	\centering
	\includegraphics[width=7.0cm]{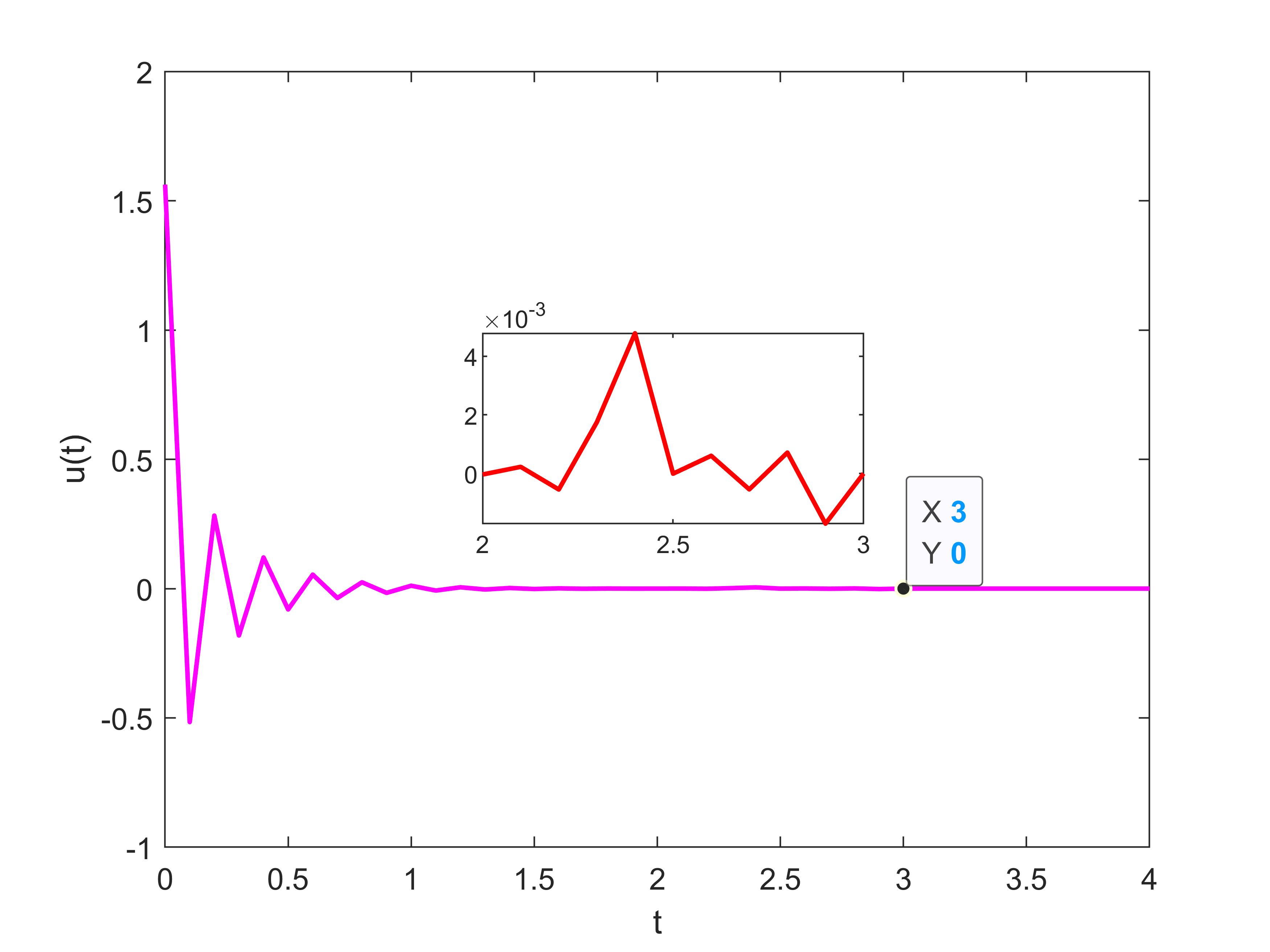}
	\caption{The evolution of the boundary controller}
	\label{fig.4}
\end{figure}\par
\indent The time evolution of the control input is shown in Fig.\ref{fig.4}, which shows that the control input is bounded and converges to zero within the prescribed time.
\newpage
\begin{figure}[!h]
	\centering
	\includegraphics[width=7.0cm]{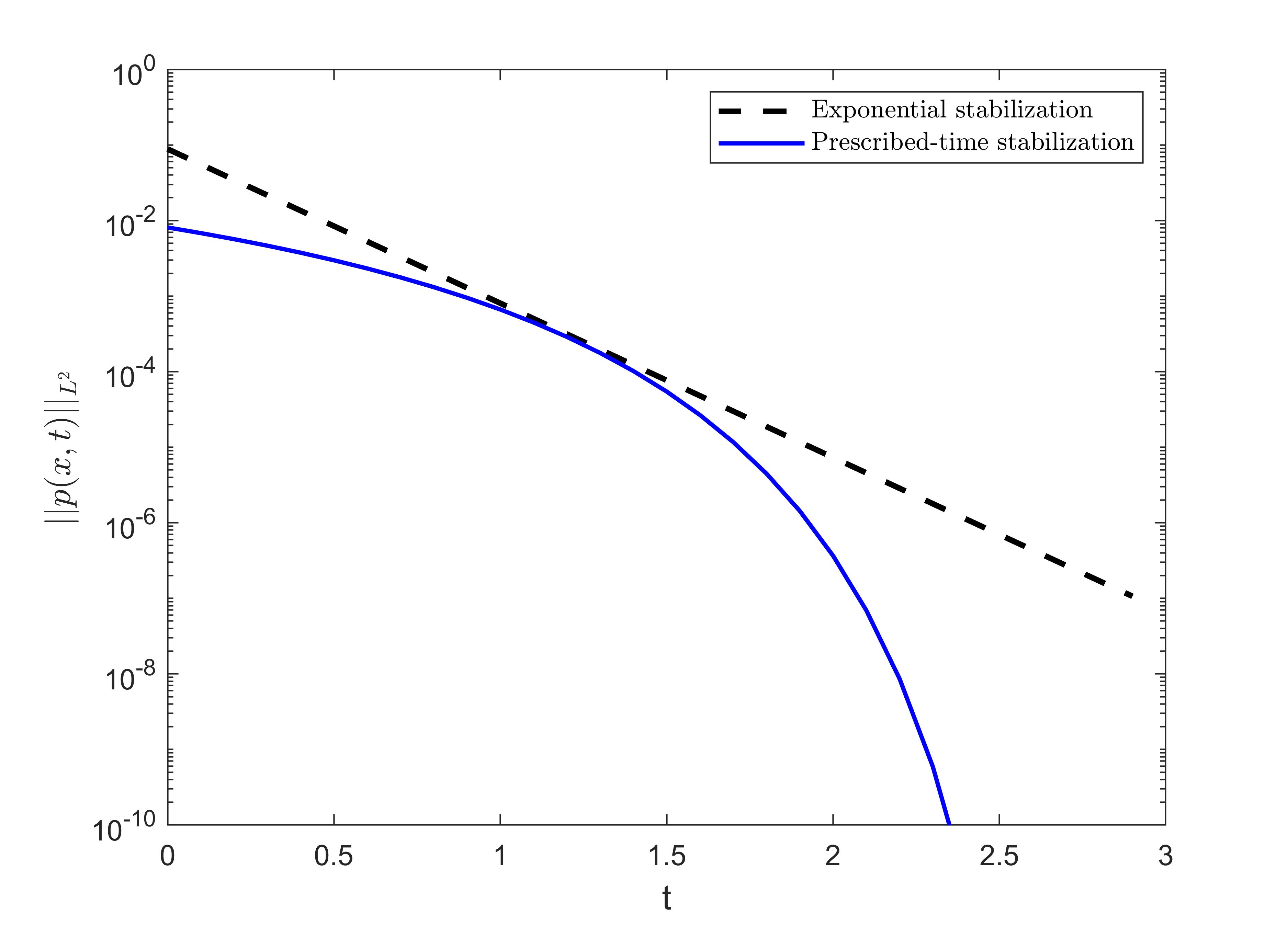}
	\caption{The trajectory of $L^{2}$-norm of $p\left ( x,t \right )$}
	\label{fig.5}
\end{figure}\par
\indent Fig.\ref{fig.5} presents the evolution of the $L^{2}$-norm of the system state under PT stabilization and exponential stabilization in a logarithmic scale. The PT controller exhibits significantly faster convergence while strictly adhering to the wave propagation constraints.
\begin{figure}[!h]
	\centering
	\includegraphics[width=7.0cm]{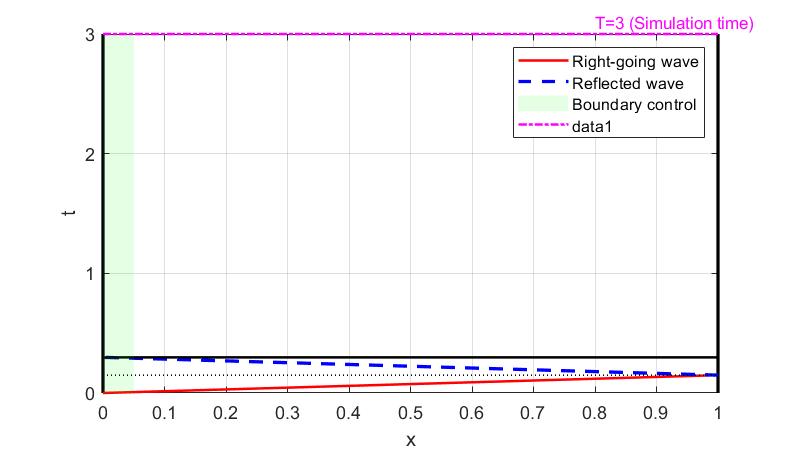}
	\caption{Wave propagation characteristics}
	\label{fig.6}
\end{figure}\par
\indent Fig.\ref{fig.6} provides a schematic of the wave propagation characteristics, control action at $x=0$ requires time $1/c$ to reach $x=1$, and the reflected wave needs another $1/c$ to return, establishing $T_{min}=2/c$. This characteristic analysis confirms our time constraint coincides with the physical wave propagation limit.
\begin{figure}[!h]
	\centering
	\includegraphics[width=7.0cm]{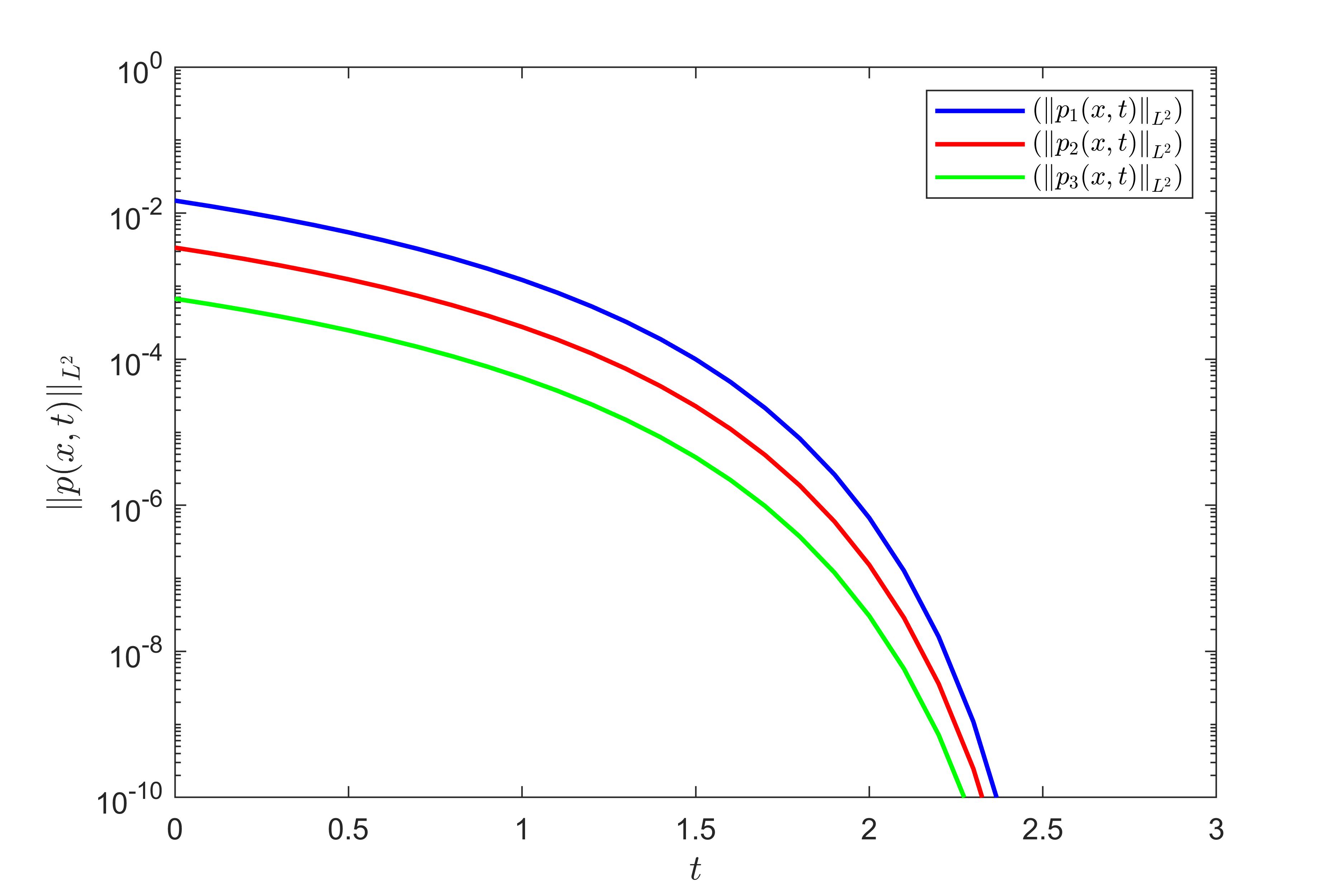}
	\caption{The trajectory of $L^{2}$-norm of $p\left ( x,t \right )$}
	\label{fig.7}
\end{figure}\par
\indent Fig.\ref{fig.7} shows that the evolution of the state of the closed-loop system under the different initial conditions converges to zero within the predefined time. The simulation result shows that the PT stabilization is independent of initial conditions of the string system.
\section{Conclusion}\label{sec6}
This paper addresses the problem of PT boundary control for flexible string systems governed by hyperbolic PDEs using a backstepping approach. A PT boundary control scheme is proposed for hyperbolic PDEs. The design incorporates a time-varying kernel function in the Volterra integral transform, enabling exact convergence within a user-specified time. The well-posedness of the kernel equations is rigorously established through successive approximation and mathematical induction. Upper bounds for the time-varying kernel functions are derived, ensuring the controller's boundedness and robustness. Simulations demonstrate the effectiveness of the proposed controller, achieving stabilization within the prescribed time while respecting the minimal controllability constraint imposed by wave propagation dynamics. Theoretical and numerical results confirm that the closed-loop system achieves PT stability without requiring distributed actuation or state observers, offering practical advantages for real-world applications. Future work will further consider the PT stabilization of perturbed flexible string systems, as well as improving noise robustness in the presence of time-varying gains.
\section*{Acknowledgements}
This work was supported by National Natural Science Foundation of China (62003189); the Natural Science Foundation of Shandong Province (Grant Nos.ZR2024QF012, ZR2022QF075, ZR2024MA060 and ZR2024QF255); the Yantai Science and Technology Innovation Development Plan (Grant 2024YT06000226).

\section*{Declarations}

\textbf{Conflict of interest} The authors declare that they have no competing interests.
\section*{CRediT authorship contribution statement}

\textbf{Chuan Zhang:} Funding acquisition; Methodology; Supervision; Writing - original draft; Writing – review \& editing. \\
\indent\textbf{He Yang:} Validation; Writing – original draft; Methodology.\\
\indent\textbf{Fei Wang:} Funding acquisition; Writing – review \& editing.\\
\indent\textbf{Tuo Zhou:}  Funding acquisition; Visualization.
\section*{Appendix A}
\textbf{Proof of Theorem 2}\\
\begin{pf}
	Make variable substitution: $\xi=x+y$, $\eta=x-y$, let $F\left ( \xi ,\eta ,t \right ) =k\left ( x,y,t \right )=k\left ( \frac{\xi +\eta }{2},\frac{\xi -\eta }{2} ,t  \right )$, then the kernel equations (\ref{shi19}) are transformed into: 
		\begin{subequations}\label{shi33}
			\begin{align}
				F_{\xi \eta } \left ( \xi ,\eta ,t \right )&=\frac{\rho _{0} }{4T_{f} }F_{tt}\left ( \xi ,\eta ,t \right ) +\frac{1}{4T_{f} }\mu \left ( t \right )F\left ( \xi ,\eta ,t \right ), \label{shi33a}  \\  
				F_{\xi} \left ( \xi ,0,t \right )&=-\frac{1}{4T_{f}}\mu \left ( t \right ),   \label{shi33b} \\
				F\left ( \xi ,\xi ,t \right )&=0,  \label{shi33c}
			\end{align}
		\end{subequations}
	where the domain of definition of $F\left ( \xi ,\eta ,t \right )$ is
	$\mathbb{D}_{1}=\left \{ \left ( \xi ,\eta ,t \right )\in \mathbb{R}^{2}\times \left [ 0,T\right ) :\eta \in \left [ 0,1 \right ],\xi \in \left [ \eta ,2-\eta  \right ] \right \} $.  \\
	\indent To solve the above equations by successive approximation, it is necessary to transform (\ref{shi33}) into an integral equation. Integrating (\ref{shi33a}) with respect to $\eta$ from 0 to $\eta $ and using (\ref{shi33b}) yields
				$\int_{0}^{\eta }F_{\xi \eta }\left ( \xi ,\eta ,t \right )d\eta=\frac{\rho _{0} }{4T_{f} }\int_{0}^{\eta }F_{tt}\left ( \xi ,s,t \right )ds+\frac{1}{4T_{f}}\int_{0}^{\eta }\mu \left ( t \right ) F\left ( \xi ,s,t \right )ds$, 
	therefore, we can get
		\begin{equation}\label{shi34}
				F_{\xi }\left ( \xi ,\eta ,t \right )=-\frac{1}{4T_{f}}\mu \left ( t \right )  +\frac{\rho _{0} }{4T_{f}}\int_{0}^{\eta }F_{tt}\left ( \xi ,s,t \right )ds+\frac{1}{4T_{f} }\int_{0}^{\eta }\mu \left ( t \right ) F\left ( \xi ,s,t \right )ds. 
		\end{equation}
	\indent Then, integrating (\ref{shi34}) with respect to $\xi$ from $\eta $ to $\xi$ and using (\ref{shi33c}) yields
		\begin{equation*}
				\int_{\eta }^{\xi } F_{\xi }\left ( \xi ,\eta ,t \right )d\xi =-\frac{1}{4T_{f} }\int_{\eta }^{\xi } \mu \left ( t \right ) d\tau+\frac{\rho _{0} }{4T_{f}}\int_{\eta }^{\xi } \int_{0}^{\eta }F_{tt}\left ( \xi ,s,t \right )ds +\frac{1}{4T_{f}}\int_{0}^{\eta }\mu \left ( t \right ) F\left ( \tau ,s,t \right )dsd\tau,
		\end{equation*}
	thereby, we obtain
		\begin{equation}\label{shi35}
				F\left ( \xi ,\eta ,t \right )=-\frac{1}{4T_{f} }\int_{\eta }^{\xi } \mu \left ( t \right ) d\tau+\frac{\rho _{0} }{4T_{f}}\int_{\eta }^{\xi } \int_{0}^{\eta }F_{tt}\left ( \xi ,s,t \right )ds  +\frac{1}{4T_{f}}\int_{\eta }^{\xi } \int_{0}^{\eta }\mu \left ( t \right ) F\left ( \tau ,s,t \right )dsd\tau. 
		\end{equation}
	\indent The solution of (\ref{shi33}) is equivalent to that of the integral equation (\ref{shi35}). Next, using successive approximation and mathematical induction to obtain the solution of the integral equation (\ref{shi35}).\\
	\indent First, take an initial condition $F^{0}\left ( \xi ,\eta ,t \right )=0$. Then, iteratively calculate equation (\ref{shi35}), the iterative equation is given as
		\begin{equation*}
				F^{n}\left ( \xi ,\eta ,t \right )=-\frac{1}{4T_{f} }\int_{\eta }^{\xi } \mu \left ( t \right ) \cdot \left ( \xi -\eta  \right )+\frac{\rho _{0} }{4T_{f}}\int_{\eta }^{\xi } \int_{0}^{\eta }F_{tt}^{n-1} \left ( \tau ,s,t \right )dsd\tau+  \frac{1}{4T_{f}}\int_{\eta }^{\xi } \int_{0}^{\eta }\mu \left ( t \right ) F^{n-1} \left ( \tau ,s,t \right )dsd\tau,
		\end{equation*}
	the continuous sequence of functions $F^{0}\left ( \xi ,\eta ,t \right )$, $F^{1}\left ( \xi ,\eta ,t \right )$, $F^{2}\left ( \xi ,\eta ,t \right )$,...,$F^{n}\left ( \xi ,\eta ,t \right )$, denoted $\left \{ F^{n} \left ( \xi ,\eta ,t \right )  \right \}$, where $n\ge 2$. If the sequence of functions converges, let $F\left ( \xi ,\eta ,t \right )=\lim\limits_{n \to \infty }F^{n}\left ( \xi ,\eta ,t \right )$, $n\ge 2$, $n\in \mathbb{N}$.\\ 
	\indent The difference between the two consecutive terms above is: $\Delta F^{n} \left ( \xi ,\eta ,t \right )=F ^{n} \left ( \xi ,\eta ,t \right )-F ^{n-1} \left ( \xi ,\eta ,t \right )$.\\  
	\indent To prove the uniform convergence of the function sequence $\left \{ F^{n}\left ( \xi ,\eta ,t \right )\right \} $, we introduce the infinite series: 
	\begin{equation}\label{shi36}
		\sum_{n=1}^{\infty } \Delta F^{n}\left ( \xi ,\eta ,t \right ),
	\end{equation}
	where it can be seen that
	\begin{align}
	\Delta F^{1}\left ( \xi ,\eta ,t \right )&=-\frac{1}{4T_{f} }\mu \left ( t \right )\cdot\left ( \xi -\eta  \right ),  \label{shi37}\\
	\Delta F^{n}\left ( \xi ,\eta ,t \right )&=\frac{\rho _{0} }{4T_{f} }\int_{\eta }^{\xi}\int_{0}^{\eta }\Delta F_{tt}^{n-1} dsd\tau\notag+\frac{1}{4T_{f} }\int_{\eta }^{\xi }\int_{0}^{\eta }\mu \left ( t \right )\notag\\
	&\quad\cdot\Delta F^{n-1}\left ( \tau ,s,t \right )dsd\tau, \quad n\ge 2.\label{shi38} 
\end{align}
	\indent It is sufficient to show that the series (\ref{shi36}) converges to obtain that the sequence $\left \{ F^{n}\left ( \xi ,\eta ,t \right ) \right \} $ is convergent. The following proves that the infinite series (\ref{shi36}) converges absolutely uniformly. \\
	\indent When $n=1$, according to (\ref{shi37}), for $\forall \left ( \xi ,\eta ,t \right )\in \mathbb{D} _{1} $, the $l\in \mathbb{N}$-order partial derivative of $\Delta F^{1} \left ( \xi ,\eta ,t \right )$ with respect to time $t$ is estimated as
		\begin{equation}\label{shi39}
				\lvert \frac{\partial^{l} }{\partial t^{l} } \Delta F^{1}\left ( \xi ,\eta ,t \right ) \rvert \le \frac{1}{4T_{f}}\int_{\eta }^{\xi }\lvert  \frac{d^{l} }{dt^{l} } \mu \left ( t \right )  \rvert d\tau \le\frac{1}{4T_{f} }\mu ^{\frac{l}{2}+1 } \left ( t \right )\left ( l+1 \right ) !\left ( \frac{1}{\mu _{0} T}  \right )^{l}\left ( \xi -\eta  \right ).  
		\end{equation}
	\indent When $n=2$, $n=3$ and $n=4$, take the $l\in \mathbb{N}$-order partial derivatives for the time $t$ for both sides of (\ref{shi38}), and use (\ref{shi39}), Lemma \ref{L1} and Lemma \ref{L2}, for $\left ( \xi ,\eta ,t \right )\in \mathbb{D}_{1}$, the estimates of the $l\in \mathbb{N}$-order partial derivatives of $\Delta F^{2} \left ( \xi ,\eta ,t \right )$, $\Delta F^{3} \left ( \xi ,\eta ,t \right )$ and $\Delta F^{4} \left ( \xi ,\eta ,t \right )$ for the time $t$ can be obtained respectively.  \\
	\indent When $n=2$, we can obtain
		\begin{equation*}
				\lvert \frac{\partial^{l} }{\partial t^{l} } \Delta F^{2}\left ( \xi ,\eta ,t \right )\rvert
				\le \left ( \frac{1}{4T_{f} } \right ) ^{2}\left ( \frac{1}{\mu _{0}T }\right )^{l+2}\mu ^{\frac{l+2}{2}+1 }  \left ( t \right )\frac{\left ( l+3 \right )! }{2\times 3}\cdot \left ( 2\times 3+\mu _{0} ^{2} T^{2}\right ) \frac{\xi \eta\left ( \xi -\eta  \right )  }{2}.
		\end{equation*}
	\indent When $n=3$, it is easy to get
		\begin{equation*}
			\begin{aligned}
				\lvert\frac{\partial^{l} }{\partial t^{l} } \Delta F^{3}\left ( \xi ,\eta ,t \right )\rvert
				\le \left ( \frac{1}{4T_{f} }\right )^{3}\left ( \left ( \frac{1}{\mu _{0}T}\right )^{l+4} \mu ^{\frac{l+4}{2}+1 }\left ( t \right ) \frac{\left ( l+5 \right )! }{4\times 5} \cdot\left ( 2\times 3 +\mu _{0}^{2} T^{2}  \right )\left ( \frac{4\times 5}{2\times 3} +\mu _{0} ^{2} T^{2}\right )\right )\frac{\left ( \xi \eta  \right ) ^{2}\left ( \xi -\eta  \right )  }{2!3!}.  
			\end{aligned}
		\end{equation*}
	\indent When $n=4$, we have access to the following equation
		\begin{equation*}
			\begin{aligned}
			\lvert\frac{\partial ^{l} }{\partial t^{l} } \Delta F^{4}\left ( \xi ,\eta ,t \right )  \rvert	&\le \left ( \frac{1}{4T_{f} }\right )^{4} \left ( \frac{1}{\mu _{0}T }\right )^{l+6}\mu ^{\frac{l+6}{2}+1 }\left ( t \right )\frac{\left ( l+7 \right )! }{6\times 7}\left ( 2\times 3+\mu _{0}^{2}T^{2}\right )\\
				&\quad \cdot\left ( \frac{4\times 5}{2\times 3}+\mu _{0}^{2}T^{2}\right ) \left ( \frac{6\times 7}{4\times 5}+\mu _{0} ^{2} T^{2}\right ) \frac{\left ( \xi \eta  \right )^{3} \left ( \xi -\eta  \right )}{3!4!},  		
			\end{aligned}
		\end{equation*}
	where $\mu \left ( t \right )=\frac{\mu _{0}^{2}T^{2} }{\left ( T-t \right )^{2} }$. For $\left ( \xi,\eta, t \right ) \in \mathbb{D}_{1}  $, the recursive form of the $l\in\mathbb{N}$-order partial derivative of $\Delta F\left ( \xi ,\eta ,t \right )$ with respect to time $t$ can be obtained by iterative computational induction as follows:
		\begin{equation}\label{shi40}
			\begin{aligned}
				\lvert\frac{\partial^{l}}{\partial t^{l}}\Delta F^{n} \left(\xi,\eta,t\right) \rvert
				&\le \left ( \frac{1}{4T_{f} }  \right ) ^{n}\left ( \frac{1}{\mu _{0}T }\right)^{l+2n-2}\mu ^{\frac{l+2n-2}{2}+1}\left ( t \right )\cdot\frac{\left ( l+2n-1 \right ) !}{\left(2n-2\right )\left (2n-1\right)}\cdot\prod_{m=1}^{n}\left[\frac{2m\left(2m+1\right)}{m\left(m+1\right)}\right.\\
				&\left.\quad+\mu _{0}^{2} T^{2}+\delta _{m,1}\left ( 3 \right )\right ]\frac{\left(\xi\eta \right )^{n-1}\left(\xi-\eta\right ) }{\left ( n-1 \right )!n! }, 
			\end{aligned}
		\end{equation}
	where $n \ge 2$, $\delta _{ij}=\left\{\begin{matrix}1,i=j \\ 0,i\ne j\end{matrix}\right.$.  \\
	\indent The following mathematical induction shows that (\ref{shi40}) holds for all positive integers $n\ge 2$.  \\
	\indent For $n=2$, it can be proved that (\ref{shi40}) holds by (\ref{shi39}). Assuming that the recursive form (\ref{shi40}) holds for positive integers $n$, the following proof that (\ref{shi40}) is suitable for positive integers $n+1$. Taking $l\in \mathbb{N} $-order partial derivatives on both sides of (\ref{shi38}) with respect to $t$, and combining (\ref{shi40}), Lemma \ref{L1} and Lemma \ref{L2} yields
		\begin{equation*}
			\begin{aligned}
				\lvert\frac{\partial ^{l} }{\partial t^{l}}\Delta F^{n+1}\left (\xi,\eta,t\right )\rvert&\le\left(\frac{1}{4T_{f}}\right)^{n+1}\left(\frac{1}{\mu _{0}T}\right )^{l+2n}\mu ^{\frac{l+2n}{2}+1}\left(t\right)\frac{\left(l+2n+1\right)!}{2n\left(2n+1\right)}   \\
				&\quad\cdot\prod_{m=1}^{n+1}\left[\frac{2m\left(2m+1\right)}{m\left(m+1\right)}+\mu _{0}^{2} T^{2} +\delta _{m,1}\left(3\right)\right]\frac{\left(\xi\eta \right )^{n}\left(\xi-\eta\right)}{n!\left ( n+1 \right )! }.
			\end{aligned}
		\end{equation*}
	\indent By mathematical induction, the recursive form (\ref{shi39}) applies for all $n\ge 2, n\in \mathbb{N}$. When $l=0$, the right side of inequality (\ref{shi39}) is a general term for the following positive term convergent series:
		\begin{equation}\label{shi41}
				\sum_{n=0}^{\infty}\left(\frac{1}{4T_{f} }\right )^{n+1}\left ( \frac{1}{\mu _{0}T }\right)^{2n}\mu ^{n+1}\left ( t \right )\cdot 2\left ( n-1 \right )!\cdot\prod_{m=1}^{n+1}\left[\frac{2m\left(2m+1\right)}{m\left(m+1\right)}+\mu _{0}^{2}T^{2} +\delta _{m,1}\left (3\right )\right ]
				\frac{\left(\xi\eta \right )^{n}\left ( \xi -\eta  \right ) }{n!\left ( n+1 \right )! }. 
		\end{equation}
	\normalsize\indent Furthermore, using the ratio convergence method for positive term series to find the radius of convergence of the power series (\ref{shi41}), we get
		\begin{equation*}
			\lim_{n\to\infty }\frac{2n!\displaystyle \prod_{m=1}^{n+2}\left[\frac{2m\left(2m+1\right)}{m\left(m+1\right)}+\mu _{0}^{2}T^{2} +\delta _{m,1}\left (3\right )\right ]\frac{1}{\left ( n+2 \right )!\left ( n+1 \right ) ! } }{2\left ( n-1 \right )!\displaystyle \prod_{m=1}^{n+1}\left[\frac{2m\left(2m+1\right)}{m\left(m+1\right)}+\mu _{0}^{2}T^{2} +\delta _{m,1}\left (3\right )\right ]\frac{1}{\left ( n+1\right )!n! }}=0. 
		\end{equation*}
	\normalsize\indent From the above analysis, it is clear that the series (\ref{shi41}) is absolutely convergent on the domain of definition $\left [ \eta ,\eta -2 \right ]\times \left [ 0,1 \right ]\times \left [ 0,T^{*}  \right ]$, where $\left [ 0,T^{*}  \right ]$, $0<T^{*}<T$ denotes any compact subset of $\left [ 0,T\right ) $. So it follows that the series (\ref{shi36}) converges absolutely uniform convergence on the domain of definition $\left [ \eta ,\eta -2 \right ]\times \left [ 0,1 \right ]\times \left [ 0,T^{*}  \right ]$. It is shown that the infinite sequence $\left \{ F^{n} \left ( \xi ,\eta ,t \right )  \right \} $ is uniformly convergent on the domain of definition $\left [ \eta ,\eta -2 \right ]\times \left [ 0,1 \right ]\times \left [ 0,T^{*}  \right ]$. Therefore, there is a unique solution $ F\left ( \xi ,\eta ,t \right )$ of the integral equation (\ref{shi35}) on $\left [ \eta ,\eta -2 \right ]\times \left [ 0,1 \right ]\times \left [ 0,T^{*}  \right ]$.\\
	\indent From the above analysis, it is clear that the equation (\ref{shi33}) has a unique solution on the domain of definition $\mathbb{D}_{1}$, and $ F\left ( \xi ,\eta ,t \right )$ is second-order continuously differentiable on $\mathbb{D}_{1}$ with respect to the spaces $\xi $ and $\eta$, as well as infinitely continuously differentiable with respect to the time $t$. This indicates that there is a unique solution $k\left ( x,y,t \right )$ of the kernel equations (\ref{shi19}) on $\mathbb{D}$, and $k\left ( x,y,t \right ) $ is second-order continuously differentiable on $\mathbb{D}$ with respect to the spaces $x$, $y$, as well as infinitely continuously differentiable with respect to the time $t$.  \\
	\indent Next, an estimate of the upper bound on the kernel functions is obtained using the series (\ref{shi41}) and the first class of first-order modified  Bessel function.  \\
	\indent Since the upper bound estimate of the obtained function $F\left ( \xi ,\eta ,t \right )$ is relative to the time-dependent sequence, the explicit expression for the upper bound estimate of $F\left ( \xi ,\eta ,t \right )$ is further given. According to
	$\lvert\displaystyle\sum_{n=0}^{\infty}\Delta F\left(\xi,\eta,t\right)\rvert\le\sum_{n=0}^{\infty}\lvert \Delta F^{n}\left ( \xi ,\eta ,t \right )\rvert$, combining with the series (\ref{shi36}) and (\ref{shi41}), the following series is obtained:
		\begin{equation}\label{shi42}
			\begin{aligned}
				\lvert\sum_{n=0}^{\infty}\Delta F^{n}\left(\xi,\eta,t\right)\rvert&\le\sum_{n=0}^{\infty}\left(\frac{1}{4T_{f}}\right)^{n+1}\mu^{n+1}\left(t\right)\left(\frac{1}{\mu _{0}T}\right)^{2n} 2\left ( n-1 \right ) !\\&\quad\cdot\prod_{m=1}^{n+1}\left [ \frac{2m\left ( 2m+1 \right)}{m\left(m+1\right)}+\mu_{0} ^{2}T^{2} +\delta_{m,1}\left(3\right)\right]\frac{\left(\xi\eta\right)^{n}\left(\xi-\eta\right)}{n!\left(n+1\right)!}.
			\end{aligned}
		\end{equation}
	To obtain an explicit expression for the upper bound of the function $ F\left ( \xi ,\eta ,t \right )$ with respect to the time-growth estimate, using the inequality deflation to further simplify (\ref{shi42}).\\
	\indent When $m=1$, it is obvious that
		\begin{equation*}
			\displaystyle \prod_{m=1}^{n+1}\left [ \frac{2m\left ( 2m+1 \right)}{m\left(m+1\right)}+\mu_{0} ^{2}T^{2} +\delta _{m,1}\left(3\right)\right]=6+\mu _{0}^{2}T^{2},
		\end{equation*}
	when $m\ge 2$, it can be deduced that
			$\frac{2m\left(2m+1\right)}{m\left(m+1\right)}=\frac{4m^{2}+2m}{m^{2}+m}=4-\frac{2m}{m^{2}+m}=4-\frac{2}{m+1}\le 4$,
	so we can get
		\begin{equation}\label{shi43}
				\prod_{m=1}^{n+1}\left[      \frac{2m\left(2m+1\right)}{m\left(m+1\right)}+\mu_{0}^{2}T^{2}+\delta_{m,1}\left(3\right)\right]\le\prod_{m=1}^{n+1}\left(6+\mu_{0}^{2}T^{2}\right)\le\left(6+\mu_{0}^{2}T^{2}\right)^{n+1}.
		\end{equation}
	\indent Therefore, for $\forall n\ge 2,n\in \mathbb{N}$, there is
		\begin{equation}\label{shi44}
				\lvert\sum_{n=0}^{\infty}\Delta F^{n}\left(\xi,\eta,t\right)\rvert
				\le\sum_{n=0}^{\infty}\left(\frac{1}{4T_{f}}\right)^{n+1}\left(\frac{1}{\mu_{0}T}\right)^{2n}\mu^{n+1}\left(t\right)\cdot 2e^{n}\cdot\left(6+\mu_{0}^{2}T^{2}\right)^{n+1}\frac{\left(\xi\eta\right)^{n}\cdot \left(\xi-\eta\right)}{n!\left(n+1\right)! }.           
		\end{equation}
	\indent  For $\xi \in \left [ \eta ,2-\eta  \right ]$, $\eta \in \left [ 0 ,1\right ] $ and $t\in\left [ 0,T\right )$, since $F\left ( \xi ,\eta ,t \right ) = \displaystyle\lim_{n\to\infty }F^{n}\left ( \xi ,\eta ,t \right ) $, $n\ge 2,n\in \mathbb{N}$ and series (\ref{shi36}), we have $\lvert  F\left ( \xi ,\eta ,t \right )\rvert= \displaystyle	\lvert\sum_{n=0}^{\infty }  \Delta F^{n}\left ( \xi ,\eta ,t \right ) \rvert$,
	and combining with (\ref{shi44}) we get
		\begin{equation}\label{shi45}
				\lvert F\left(\xi,\eta,t\right)\rvert\le\sum_{n=0}^{\infty}\left(\frac{1}{4T_{f} }\right)^{n+1}\left(\frac{1}{\mu _{0}T}\right)^{2n}\mu ^{n+1}\left(t\right)\cdot2e^{n}\cdot  \left ( 6+\mu _{0}^{2} T^{2}\right)^{n+1}  \frac{\left(\xi\eta\right)^{n}\left(\xi-\eta \right)}{n!\left(n+1\right)!}.    
		\end{equation}
	\indent To further obtain the explicit expression for $F\left ( \xi ,\eta ,t \right )$ with respect to time-varying growth, the first class of $\lambda$-order modified Bessel function is introduced:
		\begin{equation}\label{shi46}
			I_{\lambda }\left ( f \right )=\sum_{n=0}^{\infty}\frac{1}{n!\Gamma\left(n+\lambda +1\right)}\left(\frac{f}{2}\right )^{\lambda +2n}.
		\end{equation}
	\indent For the given (\ref{shi45}), we make simplified deformation: $\displaystyle a=\frac{1}{4T_{f} }\mu \left ( t \right )$, $\displaystyle b=\frac{e}{\mu _{0}^{2} T^{2}}$, $c=6+\mu _{0}^{2}T^{2}$,
	then the above equation is transformed into:
		\begin{equation}\label{shi47}
			\lvert F\left(\xi,\eta,t\right)\rvert\le \sum_{n=0 }^{\infty}\frac{a^{n+1}\cdot b^{n} \cdot c^{n+1}}{n!\left( n+1\right)!}\left ( \xi \eta  \right )^{n} \cdot \left ( \xi -\eta  \right ). 
		\end{equation}
	\indent Let $z=2\sqrt{abc\xi \eta }$, 
then (\ref{shi47}) is equivalent to
	\begin{equation}\label{shi48}
		\lvert F\left(\xi,\eta,t\right)\rvert\le\left(a\cdot c\right)^{\frac{1}{2}}\left(b\cdot\xi\eta\right)^{-\frac{1}{2}}\sum_{n=0}^{\infty}\frac{1}{n!\left(n+1\right)!}\left(\frac{z}{2}\right)^{2n+1}\left ( \xi -\eta\right ),
\end{equation}
let $\lambda=1$ in (\ref{shi46}), 
then substituting $a$, $b$, $c$ and  $z$ into (\ref{shi48}), (\ref{shi48}) can be further simplified as
\begin{equation}\label{shi49}
	\begin{aligned}
		\lvert       
		F\left(\xi,\eta,t\right)\rvert
		&\le\left ( \frac{e\xi \eta \mu \left ( t \right )\left ( 6+\mu _{0}^{2} T^{2}\right )}{4T_{f}\mu _{0} ^{2}T^{2}}\right )^{-\frac{1}{2} }\left ( \frac{\mu \left ( t \right )\left ( 6+\mu _{0}^{2} T^{2}\right)}{4T_{f} }\right )\cdot\left ( \xi -\eta  \right )\cdot I_{1}\left ( z \right ),        
	\end{aligned}
\end{equation}
where $\left ( \xi,\eta,t \right)\in \mathbb{D}_{1}$. \\
\indent Finally, by substituting the transformations $\xi =x+y$ and $\eta =x-y$ into (\ref{shi45}) and (\ref{shi49}), an estimate of the upper bound of the time-varying kernel function $k\left ( x,y,t \right ) $ can be obtained:
\begin{equation*}
	\begin{aligned}
		\lvert k\left ( x,y,t \right)\rvert&\le \sum_{n=0}^{\infty}\left(\frac{1}{4T_{f} }\right)^{n+1}\mu ^{n+1}\left ( t \right )\left(\frac{1}{\mu _{0}T }\right )^{2n}\cdot2\left ( n-1 \right )!\cdot \left( 6+\mu _{0} ^{2} T^{2}  \right )^{n+1}  \frac{\left ( \xi \eta \right)^{n} \cdot \left ( \xi -\eta  \right ) }{n!\left ( n+1 \right )! }  \\
		&\le \frac{y\mu \left ( t \right )\left ( 6+\mu _{0}^{2} T^{2}\right)}{T_{f} } \left ( \frac{e\mu \left ( t \right )\left ( 6+\mu _{0}^{2} T^{2}\right )\left ( x^{2} -y^{2}  \right ) }{T_{f}\mu _{0} ^{2}T^{2}}\right )^{-\frac{1}{2} }I_{1}\left ( z \right ),
	\end{aligned}
\end{equation*}
where $\left ( x,y,t \right )\in \mathbb{D}$. 
\end{pf}



\bibliographystyle{unsrt}
\bibliography{cas-sc-sample}

\end{document}